\documentclass[11pt,reqno]{amsart}
\usepackage{amsmath,amssymb,amsfonts,mathtools,mathrsfs,bm}
\usepackage[width=150mm,height=220mm,centering]{geometry}
\usepackage{array,tabularx,booktabs,longtable,enumitem}
\usepackage[noadjust]{cite}
\usepackage[hidelinks]{hyperref}
\allowdisplaybreaks[3]
\numberwithin{equation}{section}
\newtheorem{theorem}{Theorem}[section]
\newtheorem{lemma}[theorem]{Lemma}
\newtheorem{proposition}[theorem]{Proposition}
\newtheorem{corollary}[theorem]{Corollary}

\theoremstyle{definition}

\newtheorem{remark}[theorem]{Remark}
\newcommand{\R}{\mathbb R}
\newcommand{\dd}{\,\mathrm d}
\newcommand{\HS}{\mathrm{HS}}
\newcommand{\op}{\mathrm{op}}
\newcommand{\Tg}{T_G}
\newcommand{\one}{\mathbf 1}

\hypersetup{pdftitle={Uniqueness of Extremizers and Convergence of the Euler--Lagrange Iteration for the Critical Hardy--Littlewood--Sobolev Inequality},pdfauthor={Hanwen Lu},pdfkeywords={Hardy--Littlewood--Sobolev inequality, extremizers, total positivity, Euler--Lagrange iteration}}
\title[Uniqueness and iteration for the critical HLS inequality]
{Uniqueness of Extremizers and Convergence of the Euler--Lagrange Iteration for the Critical Hardy--Littlewood--Sobolev Inequality}
\author[H. Lu]{Hanwen Lu}
\address{School of Mathematical Sciences and LPMC,
Nankai University,
Tianjin,
China}
\email{2013537@mail.nankai.edu.cn}
\date{}
\begin{document}
\begin{abstract}
Let \(n\ge3\), \(2\le\alpha<n\), and \(1<p,r<\infty\) satisfy
\(p^{-1}+r^{-1}=1+\alpha/n\).
We prove uniqueness, up to common translations and dilations, of the normalized
nonnegative extremizing pairs for the corresponding Hardy--Littlewood--Sobolev
inequality. More generally, all positive finite-energy solutions of the
Euler--Lagrange system belong to a single symmetry orbit. The proof compares
two solutions at different relative scales. Total positivity of order two
identifies the second eigenspace of an associated compact self-adjoint
operator and forces each nonzero difference to change sign exactly once.
The opposite orientations obtained at the two ends of the scale parameter
then yield uniqueness. The same method proves radial nondegeneracy modulo
dilation. We also establish strong convergence of the normalized alternating
Euler--Lagrange iteration from every nonnegative, radially nonincreasing
initial datum to a fixed extremizing pair, without rescaling the iterates.
Square-root mass variables reduce the exclusion of several concentration
profiles to a static spectral estimate and an energy-increment argument.
The analysis includes sublinear exponents and the logarithmic borderline
decay regimes.
\end{abstract}
\maketitle

\section{Introduction and main results}\label{sec:main}

The Hardy--Littlewood--Sobolev (HLS) inequality controls a nonlocal interaction
by Lebesgue norms. Its sharp form connects fractional integration with
critical Sobolev embeddings, while its Euler--Lagrange equations describe
the associated extremal profiles. Three questions arise naturally: whether
these profiles are unique modulo the symmetries of the problem, whether the
linearized equations have additional zero modes, and whether an iteration
defined by the variational problem actually converges to an extremizer.
These questions are linked by the same obstruction: translations and
dilations preserve the functional but prevent ordinary compactness.

Let $n\ge3$, $0<\alpha<n$, $\lambda=n-\alpha$, and let $p,r>1$ satisfy
$p^{-1}+r^{-1}=1+\alpha/n$. We consider
\begin{equation}\label{eq:functional}
 \mathcal J(f,g)=
 \iint_{\R^n\times\R^n}\frac{f(x)g(y)}{|x-y|^\lambda}\dd x\dd y,
 \qquad
 C_\sharp=\sup_{\substack{f,g\ge0\\\|f\|_p=\|g\|_r=1}}\mathcal J(f,g).
\end{equation}
The HLS inequality gives $0<C_\sharp<\infty$. Lieb
\cite[Theorem~2.3(i)]{Lieb} proved attainment throughout this range and
determined the sharp constant and all extremizers in the conformally
invariant case $p=r=2n/(n+\alpha)$ \cite[Theorem~3.1]{Lieb}.
Carlen and Loss \cite{CarlenLoss1990} obtained optimizing sequences through
competing symmetries, and Frank and Lieb \cite{FrankLieb2011} subsequently
gave a rearrangement-free proof of the sharp diagonal inequality.
For unequal exponents, the two extremal components need not coincide,
and the diagonal explicit profile does not supply their classification.
Throughout this paper the Riesz kernel has exactly the normalization
in \eqref{eq:functional}.

\subsection{Classification and linearized equations}
The integral-equation approach to symmetry was developed by Chen, Li and
Ou \cite{ChenLiOu2005,ChenLiOu2006}. Their scalar classification identifies
the positive regular solutions of the conformally invariant integral
equation. For the coupled system, Hang \cite[Theorem~1.1]{Hang}
established smoothness and common-center radial symmetry of positive
finite-energy solutions, including nonlinear powers below one.
This result reduces the finite-energy classification problem to radial
profiles, but leaves the comparison of different profiles to be addressed.
Cheng, Huang and Li \cite[Section~1]{ChengHuangLi2016} formulated the
essential-uniqueness question explicitly: whether all positive critical
solutions are generated from a single pair by translations and dilations.
Our classification addresses this question in the finite-energy integral
class for $2\le\alpha<n$.
Lei and Li \cite[Theorem~1.2 and Propositions~4.2--4.5]{LeiLi}
identified their precise fast, slow, and logarithmic borderline tails
when both nonlinear powers exceed one. We use these results directly
within that range and verify the additional estimates needed below.

The Newtonian case $\alpha=2$ has a substantial independent theory.
After an amplitude normalization, its integral system becomes the critical
Lane--Emden system. Hulshof and van der Vorst
\cite[Theorems~1 and~2]{HulshofVorst1996} proved uniqueness of ground
states up to the natural symmetries and determined their asymptotic
behavior. Frank, Kim and Pistoia \cite{FrankKimPistoia2021} proved that
the ground-state linearized kernel consists precisely of the translation
and dilation modes in the corresponding energy class, and also for
decaying classical solutions. More recently, Guo, Hu, Peng and Yuan
\cite{GuoHuPengYuan2026} classified all positive classical solutions of
the critical Lane--Emden system, without an a priori finite-energy or
decay assumption. The latter result is currently a preprint.
Accordingly, our treatment of $\alpha=2$ gives a common integral-operator
proof within an established classification theory. Our aim is to carry
the comparison and radial spectral argument through the full range
$2\le\alpha<n$, including noninteger orders and sublinear powers,
and then to use it in the convergence analysis of the alternating iteration.

\subsection{Quantitative stability and concentration}
The description of extremizers also underlies quantitative remainder
estimates. Bianchi and Egnell \cite{BianchiEgnell1991} bounded the
Sobolev deficit below by the squared distance to the extremal manifold.
Recent work of Dolbeault, Esteban, Figalli, Frank and Loss
\cite{DolbeaultEtAl2025} gives constructive bounds with optimal
dimensional dependence. Chen, Lu and Tang \cite{ChenLuTang2024}
established explicit lower bounds for stability in the diagonal HLS
inequality and related fractional Sobolev inequalities.
K\"onig \cite{Konig2023,Konig2025} showed that the sharp global Sobolev
stability constant lies strictly below the local spectral threshold and
proved attainment of the corresponding stability quotient.
These results emphasize that identifying a local spectral gap and
proving a sharp global deficit inequality are distinct tasks.

A second line of work quantifies the decomposition of approximate
critical points into several bubbles. Figalli and Glaudo
\cite{FigalliGlaudo2020} obtained sharp stability results for critical
Sobolev points in low dimensions and exhibited the higher-dimensional
obstruction. Deng, Sun and Wei \cite{DengSunWei2025} established the
sharp higher-dimensional estimates for Struwe's decomposition, while
Chen, Kim and Wei \cite{ChenKimWei2025} developed the corresponding
fractional theory. The recent preprint of Chen, Lu and Tang
\cite{ChenLuTang2026} treats stability near one bubble for critical
points of the scalar HLS inequality in its natural non-Hilbertian distance
and develops a dual stability
framework. These scalar results provide the closest quantitative context
for concentration in the present problem. Here the unknowns form a coupled
system with generally unequal exponents, and the required conclusion is
convergence of an entire discrete orbit. In particular, a decomposition
of approximate stationary points must be supplemented by an argument
excluding several profiles along that orbit.

\subsection{Iteration, method, and scope of the results}
Constructive approaches to sharp inequalities provide a further motivation.
Besides the competing-symmetry scheme of Carlen and Loss
\cite{CarlenLoss1990}, Carlen, Carrillo and Loss
\cite{CarlenCarrilloLoss2010} used fast diffusion to prove the sharp
Newtonian HLS inequality. The alternating maximization studied here is
instead a nonlinear power iteration: at each half-step one maximizes over
one unit-norm component while fixing the other. Its finite-dimensional
mixed-norm predecessor is Boyd's power method \cite{Boyd}.
Gautier, Hein and Tudisco \cite{GautierHeinTudisco2021} proved global
convergence for mixed matrix norms under a projective contraction
condition. That condition is unavailable for the present convolution
kernel on the whole line: its projective cross ratios are unbounded,
and dilation generates a family of normalized fixed points.
On the whole space, monotonicity of the objective supplies neither
compactness nor convergence of the scale. The possibility of profiles
separating by arbitrarily large dilation ratios must be excluded.

Our comparison of two solutions is performed after passing to logarithmic
radial coordinates, where dilations become translations. The kernel
factorization is the logarithmic-coordinate version of the radial
Riesz-potential factorization of Lamb and Schiavone
\cite[Theorem~3.4 and Corollary~3.5]{LambSchiavone}.
Together with Karlin's composition formula
\cite[Lemma~0.1]{Karlin}, it yields order-two total positivity.
Two applications of the Perron--Jentzsch theorem
\cite[Theorem~2]{Drnovsek}, to the weighted operator and its second
exterior power, give the order-two oscillation principle
\cite[Theorem~6.2]{KushelZabreiko} in the self-adjoint setting needed here.
We verify compactness and strict positivity of the two kernels,
including at $\alpha=2$.
The resulting second-eigenfunction comparison gives uniqueness and the
radial linearized kernel. For the iteration, the
concentration--compactness strategy of Lions \cite{Lions} gives a
finite-profile decomposition. We then pass to the square roots of
the mass densities. In these Hilbert-space variables, a static
spectral estimate controls the distance to the reference profiles
by the one-step residual, and the energy gain excludes several
separating profiles. A final local argument makes the increments
of the remaining dilation parameter summable. The main conclusions are
classification modulo common symmetries and strong convergence, from
arbitrary nonnegative radially nonincreasing initial data, to one fixed
extremizing pair, without rescaling successive iterates.

We now state the results precisely. The first theorem concerns normalized
nonnegative extremizing pairs; the subsequent statements treat all
positive finite-energy solutions and the alternating iteration.

\begin{theorem}[Uniqueness of extremizing pairs]\label{thm:main}
Suppose that \(n\ge3\), \(2\le\alpha<n\), and
\begin{equation}\label{eq:critical}
 1<p,r<\infty,\qquad \frac1p+\frac1r=1+\frac{\alpha}{n}.
\end{equation}
Fix a nonnegative normalized extremizing pair \((f_*,g_*)\) with common
radial center at the origin, so that
\[
 \|f_*\|_p=\|g_*\|_r=1,\qquad \mathcal J(f_*,g_*)=C_\sharp.
\]
For every nonnegative normalized extremizing pair \((f,g)\), there exist
\(x_0\in\R^n\) and \(\rho>0\) such that, almost everywhere,
\begin{equation}\label{eq:uniqueness}
 f(x)=\rho^{n/p}f_*\bigl(\rho(x-x_0)\bigr),\qquad
 g(x)=\rho^{n/r}g_*\bigl(\rho(x-x_0)\bigr).
\end{equation}
Conversely, every pair in \eqref{eq:uniqueness} is a nonnegative normalized
extremizing pair.
\end{theorem}

Both components in \eqref{eq:uniqueness} have the same center and scale;
their norm-preserving amplitude factors differ when \(p\ne r\).
The theorem neither asserts \(f_*=g_*\) nor requires an explicit formula
for either profile. Without norm normalization, all nonzero nonnegative
extremizing pairs are
\[
 f(x)=A\rho^{n/p}f_*\bigl(\rho(x-x_0)\bigr),\qquad
 g(x)=B\rho^{n/r}g_*\bigl(\rho(x-x_0)\bigr),\quad A,B>0,
\]
where \(A=\|f\|_p\) and \(B=\|g\|_r\) are independent.
The converse follows by a change of variables, since
\(n/p+n/r-2n+\lambda=0\).
The parameters \(n,\alpha,p,r\) remain fixed throughout the comparison.

\subsection{The Euler--Lagrange system and amplitude normalization}
Set
\begin{equation}\label{eq:parameters}
 P=\frac1{p-1},\qquad Q=\frac1{r-1},\qquad
 a=\frac n{P+1},\qquad b=\frac n{Q+1}.
\end{equation}
Since \(p^{-1}=P/(P+1)\) and \(r^{-1}=Q/(Q+1)\),
\eqref{eq:critical} is equivalent to
\[
 \frac1{P+1}+\frac1{Q+1}=1-\frac\alpha n,\qquad
 PQ-1=\frac\alpha n(P+1)(Q+1)>0.
\]
In particular,
\begin{equation}\label{eq:identities}
 a+b=\lambda,\qquad aP=b+\alpha,\qquad
 bQ=a+\alpha,\qquad PQ>1.
\end{equation}
The nonlinear exponents \(P,Q\) need only be positive and satisfy these
relations. Thus \(P<1\), corresponding to \(p>2\), and \(Q<1\),
corresponding to \(r>2\), are included.
We write \(w\asymp h\) for two-sided comparison by positive finite
constants, and \(w\sim h\) for \(w/h\to1\).
Unless stated otherwise, constants may depend on the fixed parameters
and the profiles under consideration, but not on the limiting variable.

Define
\[
 (Th)(x)=\int_{\R^n}|x-y|^{-\lambda}h(y)\dd y.
\]
A normalized extremizing pair satisfies
\begin{equation}\label{eq:ELnorm}
 Tg=C_\sharp f^{p-1},\qquad Tf=C_\sharp g^{r-1}.
\end{equation}
Indeed, optimizing one component and using equality in H\"older's
inequality gives these equations; integration against that component
identifies the multiplier as \(C_\sharp\).
Strict positivity of the kernel implies positivity of both components.
With
\begin{equation}\label{eq:amplitude}
 U=C_\sharp^{-(Q+1)/(PQ-1)}f^{p-1},\qquad
 V=C_\sharp^{-(P+1)/(PQ-1)}g^{r-1},
\end{equation}
the system becomes
\begin{equation}\label{eq:physicalEL}
 U=T(V^Q),\qquad V=T(U^P).
\end{equation}
Moreover,
\begin{equation}\label{eq:massphysical}
 \int_{\R^n}U^{P+1}\dd x
 =\int_{\R^n}V^{Q+1}\dd x
 =C_\sharp^{-n/\alpha},\qquad
 \frac{(P+1)(Q+1)}{PQ-1}=\frac n\alpha.
\end{equation}
Thus normalized extremizing pairs have identical power integrals.

By Hang \cite[Theorem~1.1]{Hang}, positive finite-energy solutions of
\eqref{eq:physicalEL} are smooth and strictly radially decreasing about
a common center. The hypotheses follow from
\eqref{eq:identities} and \(U\in L^{P+1}\), \(V\in L^{Q+1}\);
in particular they allow \(P<1\) or \(Q<1\).
We translate each pair so that its center is the origin.

\subsection{The radial cylindrical transformation}
Writing \(R=|x|=e^t\), define
\begin{equation}\label{eq:EF}
 u(t)=e^{at}U(e^t),\qquad v(t)=e^{bt}V(e^t),\qquad t\in\R.
\end{equation}
Let \(\omega_{n-1}=|\mathbb S^{n-1}|\). Polar coordinates give
\begin{equation}\label{eq:cylEL}
 u=\Tg(v^Q),\qquad v=\Tg^*(u^P),\qquad
 \Tg h=G*h,\quad G^*(t)=G(-t),
\end{equation}
where
\begin{equation}\label{eq:G}
 G(t)=e^{(a-b)t/2}
 \int_{\mathbb S^{n-1}}(2\cosh t-2\theta_1)^{-\lambda/2}\dd\theta.
\end{equation}
Here \(\theta=(\theta_1,\ldots,\theta_n)\in\mathbb S^{n-1}\),
\(\theta_1=\theta\cdot e_1\) is its first coordinate, and
\(\dd\theta\) denotes surface measure. Also \(\Tg^*h=G^**h\).
For example,
\[
 |e^t e_1-e^z\theta|^{-\lambda}
 =e^{-\lambda(t+z)/2}(2\cosh(t-z)-2\theta_1)^{-\lambda/2}.
\]
Substitution into the first integral equation and the identity
\(n-bQ=b\) yield the first equation in \eqref{eq:cylEL};
the second follows in the same way.
The cylindrical profiles satisfy
\begin{equation}\label{eq:cylmass}
 \int_\R u^{P+1}\dd t=\int_\R v^{Q+1}\dd t=S,\qquad
 S=\omega_{n-1}^{-1}C_\sharp^{-n/\alpha}.
\end{equation}
A common dilation in physical space corresponds precisely to
\begin{equation}\label{eq:translation}
 (u,v)(t)\longmapsto (u,v)(t+s),\qquad s\in\R.
\end{equation}
It therefore remains to prove uniqueness of extremizing profiles
satisfying \eqref{eq:cylmass}, modulo common translations.

\subsection{Classification and convergence}
The uniqueness argument applies to all positive finite-energy solutions.

\begin{theorem}[Classification of positive finite-energy solutions]\label{thm:allEL}
Fix the parameters in \eqref{eq:parameters}--\eqref{eq:identities} and
the kernel in \eqref{eq:G}. For \(i=0,1\), suppose that
\[
 u_i>0,\quad v_i>0,\qquad
 u_i\in L^{P+1}(\R),\quad v_i\in L^{Q+1}(\R),
\]
and
\[
 u_i=G*v_i^Q,\qquad v_i=G^**u_i^P.
\]
There exists \(s_*\in\R\) such that
\[
 u_1(t+s_*)=u_0(t),\qquad v_1(t+s_*)=v_0(t)\quad(t\in\R).
\]
No equality of the two pairs' power integrals, or extremizing property,
is assumed. Consequently, positive solutions of
\[
 U=T(V^Q),\quad V=T(U^P),\qquad
 U\in L^{P+1}(\R^n),\quad V\in L^{Q+1}(\R^n)
\]
are unique up to common translations and
\((U,V)(x)\mapsto(\rho^aU(\rho x),\rho^bV(\rho x))\).
\end{theorem}

For \(f_0\ge0\) with \(\|f_0\|_p=1\), define the alternating iteration by
\begin{equation}\label{eq:mainiteration}
 g_k=\frac{(Tf_k)^Q}{\|(Tf_k)^Q\|_r},\qquad
 f_{k+1}=\frac{(Tg_k)^P}{\|(Tg_k)^P\|_p},\qquad k=0,1,\ldots.
\end{equation}
For fixed \(f_k\), the first step is the unique maximizer of
\(\mathcal J(f_k,g)\) on the nonnegative \(L^r\) unit sphere.
The second step has the corresponding interpretation with \(g_k\) fixed.
Both assertions follow from equality in H\"older's inequality.

\begin{theorem}[Convergence of the alternating iteration]\label{thm:mainiteration}
Under the assumptions of Theorem~\ref{thm:main}, suppose that
\(f_0\ge0\), \(\|f_0\|_p=1\), and \(f_0\) is radially nonincreasing
about \(x_0\). Every step of \eqref{eq:mainiteration} is well defined.
There exists a fixed \(\rho_\infty\in(0,\infty)\), depending on the
initial datum, such that
\begin{align}
 &\left\|f_k-\rho_\infty^{n/p}
     f_*\bigl(\rho_\infty(\,\cdot-x_0)\bigr)\right\|_{L^p(\R^n)}
 \notag\\
 &\quad+
 \left\|g_k-\rho_\infty^{n/r}
     g_*\bigl(\rho_\infty(\,\cdot-x_0)\bigr)\right\|_{L^r(\R^n)}
 \longrightarrow0.\label{eq:mainiterationlimit}
\end{align}
Furthermore, \(\mathcal J(f_k,g_k)\) is nondecreasing and converges to
\(C_\sharp\). No smoothness, strict positivity, moment condition,
prescribed tail behavior, or proximity to the optimal value is required
of \(f_0\).
\end{theorem}

The scale in \eqref{eq:mainiterationlimit} is a single limiting number,
not a sequence of scales chosen to normalize successive iterates.
Translations introduced in the proof are only used to compare
cylindrical profiles. A nonzero unnormalized initial datum can first be
divided by its \(L^p\) norm.

Sections~\ref{sec:direct-kernel}--\ref{sec:secant} establish the kernel,
tail, and spectral estimates underlying the classification.
The classification and nondegeneracy arguments are then followed by
the iteration analysis. A separate part of that analysis excludes mass
retained in regions whose radii have unbounded ratios: monotonicity of
the objective and uniqueness of extremizers alone do not imply convergence.

\section{Factorization and order-two total positivity of the cylindrical kernel}\label{sec:direct-kernel}

Throughout this section we use the unnormalized Riesz kernel
\[
(Th)(x)=\int_{\mathbb R^n}|x-y|^{-\lambda}h(y)\,\mathrm dy.
\]
Assume that
\[
n\ge3,\qquad 2\le\alpha=n-\lambda<n,\qquad
a,b>0,\qquad a+b=\lambda,\qquad \sigma=\frac{\alpha}{2}.
\]
Write $\omega_{n-1}=|\mathbb S^{n-1}|$. The cylindrical kernel associated with the radial transformation is
\begin{equation}\label{eq:direct-kernel-G}
G(t)=e^{(a-b)t/2}
\int_{\mathbb S^{n-1}}
(2\cosh t-2\omega_1)^{-\lambda/2}\,\mathrm d\omega .
\end{equation}
The one-sided factorization below is the logarithmic-coordinate form of the radial Riesz-potential factorization into Erd\'elyi--Kober operators of Lamb and Schiavone \cite[Corollary~3.5, formula~(3.12)]{LambSchiavone}. We record its normalization and the kernel estimates needed here, and then establish strict order-two total positivity for the two-step kernel.

\begin{lemma}\label{lem:direct-kernel-factor}
For $c>0$, define
\begin{equation}\label{eq:direct-kernel-side}
h_{c,+}(t)=
\frac{e^{-ct}(1-e^{-2t})^{\sigma-1}}{\Gamma(\sigma)}
\mathbf 1_{\{t>0\}},
\qquad h_{c,-}(t)=h_{c,+}(-t).
\end{equation}
Then
\begin{equation}\label{eq:direct-kernel-factor}
G=\gamma_{n,\alpha}\,h_{a,-}*h_{b,+},
\qquad
\gamma_{n,\alpha}
=\frac{4\pi^{n/2}\Gamma(\alpha/2)}{\Gamma(\lambda/2)}.
\end{equation}
Moreover, $G$ is strictly positive, belongs to $L^1(\mathbb R)\cap L^\infty(\mathbb R)$, and satisfies
\begin{align}
G(t)&=\omega_{n-1}e^{-bt}
       \bigl(1+O(e^{-2t})\bigr), &&t\longrightarrow+\infty,\label{eq:direct-kernel-tail-plus}\\
G(t)&=\omega_{n-1}e^{at}
       \bigl(1+O(e^{2t})\bigr), &&t\longrightarrow-\infty.\label{eq:direct-kernel-tail-minus}
\end{align}
It is locally absolutely continuous, with
\begin{equation}\label{eq:direct-kernel-derivative}
-bG(t)\le G'(t)\le aG(t)
\quad\text{for almost every }t\in\mathbb R.
\end{equation}
If $\alpha>2$, then $G\in C^1(\mathbb R)$; if $\alpha=2$, then
\begin{equation}\label{eq:direct-kernel-newton}
G(t)=\omega_{n-1}
\begin{cases}
e^{at},&t\le0,\\
e^{-bt},&t\ge0.
\end{cases}
\end{equation}
\end{lemma}

\begin{proof}
Near $t=0$, the singularity in the spherical integral is $|\omega-e_1|^{-\lambda}$. Since $\lambda\le n-2<n-1$, it is integrable with respect to surface measure. The inequality $2\cosh t-2\omega_1\ge2-2\omega_1$ therefore allows dominated convergence, proving continuity and local boundedness of $G$. For $t>0$, rewrite \eqref{eq:direct-kernel-G} as
\[
G(t)=e^{-bt}\int_{\mathbb S^{n-1}}
(1-2e^{-t}\omega_1+e^{-2t})^{-\lambda/2}\,\mathrm d\omega.
\]
A uniform Taylor expansion in $e^{-t}$, together with $\int_{\mathbb S^{n-1}}\omega_1\,\mathrm d\omega=0$, gives \eqref{eq:direct-kernel-tail-plus}. Replacing $t$ by $-t$ gives \eqref{eq:direct-kernel-tail-minus}. These estimates also imply $G\in L^1\cap L^\infty$; strict positivity follows directly from its integral representation.

Use the Fourier transform convention $\widehat f(\xi)=\int_{\mathbb R}e^{-it\xi}f(t)\,\mathrm dt$. Substituting $s=e^{-2t}$ in \eqref{eq:direct-kernel-side} and applying the beta integral yields
\begin{align}
\widehat h_{c,+}(\xi)
&=\frac12\,
\frac{\Gamma((c+i\xi)/2)}
     {\Gamma((c+\alpha+i\xi)/2)},\label{eq:direct-kernel-side-fourier}\\
\widehat h_{c,-}(\xi)
&=\frac12\,
\frac{\Gamma((c-i\xi)/2)}
     {\Gamma((c+\alpha-i\xi)/2)}.\nonumber
\end{align}
On the other hand, setting $r=e^t$ in \eqref{eq:direct-kernel-G} gives
\begin{equation}\label{eq:direct-kernel-mellin}
\widehat G(\xi)=
\int_{\mathbb R^n}
|x|^{a-n-i\xi}|x-e_1|^{-\lambda}\,\mathrm dx.
\end{equation}
Here we used $|r\omega-e_1|^2=r(2\cosh t-2\omega_1)$. The integral is absolutely convergent at $0$, at $e_1$, and at infinity.

We now apply the radial Mellin multiplier formula of Lamb and Schiavone \cite[Theorem~3.4, formula~(3.11)]{LambSchiavone}. Their order parameter is $\sigma=\alpha/2$, and their normalized radial Riesz operator is $R_n^\alpha=C_{n,\alpha}T$, where
\[
C_{n,\alpha}=2^{-\alpha}\pi^{-n/2}
\frac{\Gamma(\lambda/2)}{\Gamma(\alpha/2)}.
\]
Set $s=a+\alpha-i\xi=n-b-i\xi$. The required Mellin strip $\alpha<\operatorname{Re}s<n$ holds because $a,b>0$. By the radial change of variables underlying \eqref{eq:direct-kernel-mellin}, the multiplier of $r^{-\alpha}R_n^\alpha$ at this $s$ equals $C_{n,\alpha}\widehat G(\xi)$. Dividing the cited multiplier by $C_{n,\alpha}$ therefore gives
\begin{equation}\label{eq:direct-kernel-fourier}
\widehat G(\xi)
=\frac{\pi^{n/2}\Gamma(\alpha/2)}{\Gamma(\lambda/2)}
\frac{\Gamma((a-i\xi)/2)\Gamma((b+i\xi)/2)}
     {\Gamma((a+\alpha-i\xi)/2)
      \Gamma((b+\alpha+i\xi)/2)}.
\end{equation}
Comparison with \eqref{eq:direct-kernel-side-fourier} and uniqueness of the Fourier transform on $L^1$ prove \eqref{eq:direct-kernel-factor}. Since $h_{c,\pm}\in L^1\cap L^\infty$, the convolution on the right is continuous as well, so the identity holds everywhere.

Finally, in the sense of distributions,
\begin{equation}\label{eq:direct-kernel-side-derivative}
h'_{b,+}\ge-bh_{b,+},
\qquad
h'_{a,-}\le ah_{a,-}.
\end{equation}
For $\sigma>1$, this follows by direct differentiation, and both derivatives belong to $L^1(\mathbb R)$. For $\sigma=1$, the precise identities are
\[
h'_{b,+}=\delta_0-bh_{b,+},
\qquad
h'_{a,-}=ah_{a,-}-\delta_0.
\]
Convolving each inequality in \eqref{eq:direct-kernel-side-derivative} with the other nonnegative factor yields \eqref{eq:direct-kernel-derivative}. For $\sigma>1$, convolution of an $L^1$ derivative with the bounded factor is continuous; hence $G\in C^1$. For $\sigma=1$, direct evaluation of the convolution, using $\gamma_{n,2}/(a+b)=\omega_{n-1}$, gives \eqref{eq:direct-kernel-newton}.
\end{proof}

\begin{proposition}\label{prop:direct-kernel-tp2}
For every $t_1<t_2$ and $x_1<x_2$,
\begin{equation}\label{eq:direct-kernel-tp2}
\det\bigl[G(t_i-x_j)\bigr]_{i,j=1}^2\ge0.
\end{equation}
If $\alpha>2$, the inequality in \eqref{eq:direct-kernel-tp2} is strict. If $\alpha=2$, the one-step kernel has order-two total positivity, but not strict order-two total positivity.
\end{proposition}

\begin{proof}
On $(0,\infty)$,
\[
\frac{\mathrm d^2}{\mathrm dt^2}\log h_{c,+}(t)
=-4(\sigma-1)\frac{e^{2t}}{(e^{2t}-1)^2}\le0.
\]
Extending the function by zero outside its support gives a log-concave function; its reflection $h_{c,-}$ has the same property. The classical log-concavity criterion for order-two Toeplitz kernels, recalled in Belton, Guillot, Khare and Putinar \cite[Section~8]{BeltonEtAl2023}, therefore gives order-two total positivity of $h_{c,\pm}(t-x)$. When $\sigma>1$, strict log-concavity gives a strictly positive determinant whenever all four arguments lie in the interior of the support.

The continuous Cauchy--Binet composition formula of Karlin \cite[Lemma~0.1, formula~(0.4)]{Karlin} gives
\begin{equation}\label{eq:direct-kernel-cauchy-binet}
\begin{aligned}
\det\bigl[G(t_i-x_j)\bigr]_{i,j=1}^2
={}&\gamma_{n,\alpha}^{\,2}
\int_{z_1<z_2}
\det\bigl[h_{a,-}(t_i-z_j)\bigr]_{i,j=1}^2\\
&\hspace{18mm}\cdot
\det\bigl[h_{b,+}(z_i-x_j)\bigr]_{i,j=1}^2
\,\mathrm dz_1\,\mathrm dz_2.
\end{aligned}
\end{equation}
After expanding both determinants, the absolute integral of each term is bounded by a product of finite convolution integrals; all interchanges are therefore justified. The integrand is nonnegative, proving \eqref{eq:direct-kernel-tp2}. For $\sigma>1$, both determinants are strictly positive on the open set
\[
z_2>z_1>\max\{t_2,x_2\},
\]
so the integral is strictly positive.

For $\sigma=1$, \eqref{eq:direct-kernel-newton} shows that whenever all four differences $t_i-x_j$ are positive, or all are negative, the matrix has rank one and its determinant vanishes.
\end{proof}

\begin{proposition}\label{prop:direct-kernel-gram}
Let $m:\mathbb R\to(0,\infty)$ be measurable and suppose that, for each $t\in\mathbb R$,
\begin{equation}\label{eq:direct-kernel-gram-domain}
\int_{\mathbb R}G(t-z)^2m(z)\,\mathrm dz<\infty.
\end{equation}
Define
\begin{equation}\label{eq:direct-kernel-gram}
K_m(t,x)=\int_{\mathbb R}G(t-z)m(z)G(x-z)\,\mathrm dz.
\end{equation}
For every $2\le\alpha<n$, the kernel $K_m$ has strict order-two total positivity:
\begin{equation}\label{eq:direct-kernel-gram-strict}
\det\bigl[K_m(t_i,x_j)\bigr]_{i,j=1}^2>0
\qquad(t_1<t_2,\ x_1<x_2).
\end{equation}
The conclusion remains valid after multiplying \eqref{eq:direct-kernel-gram} by arbitrary strictly positive row and column factors. If \eqref{eq:direct-kernel-gram-domain} holds only for almost every $t$, the conclusion holds on the corresponding set of full measure.
\end{proposition}

\begin{proof}
The Cauchy--Schwarz inequality first ensures finiteness of \eqref{eq:direct-kernel-gram}. Another application of Karlin's continuous Cauchy--Binet identity \cite[Lemma~0.1, formula~(0.4)]{Karlin} gives
\begin{equation}\label{eq:direct-kernel-gram-cauchy-binet}
\begin{aligned}
\det\bigl[K_m(t_i,x_j)\bigr]_{i,j=1}^2
={}&\int_{z_1<z_2}
\det\bigl[G(t_i-z_j)\bigr]_{i,j=1}^2\\
&\quad\cdot
\det\bigl[G(x_i-z_j)\bigr]_{i,j=1}^2
m(z_1)m(z_2)\,\mathrm dz_1\,\mathrm dz_2.
\end{aligned}
\end{equation}
Absolute integrability again follows by expanding the determinants and applying Cauchy--Schwarz. Proposition~\ref{prop:direct-kernel-tp2} makes the integrand nonnegative. For $\alpha>2$, both determinants are strictly positive throughout $z_1<z_2$, proving the assertion.

Now suppose $\alpha=2$ and take
\[
z_1<\min\{t_1,x_1\},
\qquad z_2>\max\{t_2,x_2\}.
\]
By \eqref{eq:direct-kernel-newton},
\[
\det\bigl[G(t_i-z_j)\bigr]_{i,j=1}^2
=\omega_{n-1}^{\,2}
e^{bz_1-az_2+at_1-bt_2}
\bigl(e^{\lambda(t_2-t_1)}-1\bigr)>0.
\]
Replacing $t_i$ by $x_i$ gives strict positivity of the other determinant. This integration region has positive measure, and $m>0$ almost everywhere, so \eqref{eq:direct-kernel-gram-cauchy-binet} is strictly positive. Positive row and column factors multiply the determinant by a positive number and do not change its sign.
\end{proof}

\section{Tail estimates and square integrability for fixed parameters}\label{sec:direct-tail}

Fix parameters
\[
 a,b,P,Q>0,\qquad
 a(P+1)=b(Q+1)=n,\qquad
 aP=b+\alpha,\quad bQ=a+\alpha,\quad PQ>1.
\]
Let positive functions $u,v$ satisfy
\begin{equation}\label{eq:direct-tail-system}
 u=G*v^Q,\qquad v=G^**u^P,\qquad
 u\in L^{P+1}(\R),\quad v\in L^{Q+1}(\R),
 \qquad G^*(t)=G(-t).
\end{equation}
The preceding kernel representation shows that $G$ is positive, bounded, and locally absolutely continuous, and that there are $k_-,k_+>0$ such that
\begin{equation}\label{eq:direct-tail-kernel}
 \begin{gathered}
 G(t)\sim k_-e^{at}\quad(t\to-\infty),\qquad
 G(t)\sim k_+e^{-bt}\quad(t\to+\infty),\\
 -bG\le G'\le aG\quad\text{almost everywhere}.
 \end{gathered}
\end{equation}
These properties hold for every $2\le\alpha<n$. Constants in this section may depend on the fixed parameters and on the solution under consideration.

\subsection{Decay estimates of Lei and Li}

For $P,Q>1$, the pointwise decay estimates used below are direct applications of Lei--Li \cite[Theorem~1.2 and Propositions~4.2--4.5]{LeiLi}. To specify the parameter correspondence, define the physical profiles by
\[
 U(R)=R^{-a}u(\log R),\qquad
 V(R)=R^{-b}v(\log R),\qquad R>0.
\]
The radial change of variables gives
\[
 U=T(V^Q),\qquad V=T(U^P),\qquad
 U\in L^{P+1}(\mathbb R^n),\quad
 V\in L^{Q+1}(\mathbb R^n).
\]
Here $T$ has the unnormalized kernel $|x|^{-\lambda}$, as in \eqref{eq:physicalEL}. Lei--Li's nonlinear exponents $(p,q)$ are $(P,Q)$, and their integrability exponents satisfy
\[
 r_0=\frac{n(PQ-1)}{\alpha(Q+1)}=P+1,\qquad
 s_0=\frac{n(PQ-1)}{\alpha(P+1)}=Q+1.
\]
Thus these finite-energy solutions belong to their integrable-solution class. In the ordering $P\le Q$, their results state that, as $R\to\infty$,
\begin{equation}\label{eq:direct-published-physical-tails}
 U(R)\sim C_U R^{-\lambda},\qquad
 V(R)\sim
 \begin{cases}
 C_V R^{-\lambda},&P\lambda>n,\\
 C_V R^{-\lambda}\log R,&P\lambda=n,\\
 C_V R^{-(P\lambda-\alpha)},&P\lambda<n,
 \end{cases}
 \qquad C_U,C_V>0.
\end{equation}
For $P>Q$, interchange $(U,P,a)$ and $(V,Q,b)$ before applying this statement. The identities
\[
 P\lambda-n=bP-a,\qquad Q\lambda-n=aQ-b
\]
translate these alternatives into the cylindrical rates $\rho_u,\rho_v$ and resonance factors stated below. This citation is used only when $P,Q>1$, as required by the cited results.

At the other end, Hang's regularity theorem \cite[Theorem~1.1]{Hang} applies for all the parameters considered here. It gives smooth positive radial extensions of $U,V$ across the origin. Consequently,
\begin{equation}\label{eq:direct-tail-origin}
 u(t)=U(0)e^{at}(1+o(1)),\qquad
 v(t)=V(0)e^{bt}(1+o(1))
 \quad(t\to-\infty).
\end{equation}

\subsection{Preliminary decay estimates}

\begin{lemma}\label{lem:direct-tail-preliminary}
The solutions $u,v$ are strictly positive $C^1$ functions, vanish at both ends of the real line, and satisfy
\begin{equation}\label{eq:direct-tail-logbound}
 -b\le \frac{u'}u\le a,\qquad
 -a\le \frac{v'}v\le b.
\end{equation}
Set
\[
 \rho_u=\min\{b,aQ\},\qquad
 \rho_v=\min\{a,bP\}.
\]
For every sufficiently small $\varepsilon>0$, there exists $C_\varepsilon<\infty$ such that
\begin{align}
 u(t)+|u'(t)|&\le C_\varepsilon e^{-(\rho_u-\varepsilon)t},
 &
 v(t)+|v'(t)|&\le C_\varepsilon e^{-(\rho_v-\varepsilon)t}
 &&(t\ge0),\label{eq:direct-tail-pre-right}\\
 u(t)+|u'(t)|&\le C_\varepsilon e^{(a-\varepsilon)t},
 &
 v(t)+|v'(t)|&\le C_\varepsilon e^{(b-\varepsilon)t}
 &&(t\le0).\label{eq:direct-tail-pre-left}
\end{align}
\end{lemma}

\begin{proof}
By \eqref{eq:direct-tail-kernel}, $G,G'\in L^s(\R)$ for every $1\le s<\infty$. Since $v^Q\in L^{(Q+1)/Q}$ and $u^P\in L^{(P+1)/P}$, H\"older's inequality, continuity of translations in $L^s$, and compactly supported approximation give $u,v\in C^1(\R)\cap C_0(\R)$, with
\[
 u'=G'*v^Q,\qquad v'=(G^*)'*u^P.
\]
Positivity of $G$ gives $u,v>0$. Convolving the inequalities in \eqref{eq:direct-tail-kernel} with the nonnegative sources gives \eqref{eq:direct-tail-logbound}.

The left-hand rates in \eqref{eq:direct-tail-pre-left} follow from \eqref{eq:direct-tail-origin} and \eqref{eq:direct-tail-logbound}. When $P,Q>1$, the right-hand rates in \eqref{eq:direct-tail-pre-right} follow directly from \eqref{eq:direct-published-physical-tails}, with the components interchanged when necessary; a small loss in the exponent absorbs the logarithmic factor. The derivative bounds then follow from \eqref{eq:direct-tail-logbound}.

It remains to prove the right-hand rates when $\min\{P,Q\}\le1$. We first obtain decay at some positive exponential rate. Define
\[
 M_u(R)=\sup_{|t|\ge R}u(t),\qquad
 M_v(R)=\sup_{|t|\ge R}v(t).
\]
Both are nonincreasing in $R$ and tend to zero. Fix $\theta\in(0,1)$ and split the first convolution into $|z|\ge\theta R$ and $|z|<\theta R$. For $|t|\ge R$, the first part is bounded by $\|G\|_1M_v(\theta R)^Q$. In the second part, $|t-z|\ge(1-\theta)R$, so it is bounded by
\[
 \|G\mathbf1_{\{|\,\cdot\,|\ge(1-\theta)R\}}\|_{Q+1}
 \|v^Q\|_{(Q+1)/Q}\le Ce^{-cR}.
\]
The same decomposition in the second equation gives
\[
 M_u(R)\le C M_v(\theta R)^Q+Ce^{-cR},\qquad
 M_v(R)\le C M_u(\theta R)^P+Ce^{-cR}.
\]
Insert the second inequality, evaluated at $\theta R$, into the first. For any fixed $\gamma>0$, the inequality $(x+y)^\gamma\le \max\{1,2^{\gamma-1}\}(x^\gamma+y^\gamma)$ justifies this substitution even for exponents below one. We obtain
\[
 M_u(R)\le A M_u(\theta^2R)^{PQ}+Ae^{-\kappa R}.
\]
Choose $\theta$ so that $\theta^2PQ>1$. The elementary iteration below gives $M_u(R)\le Ce^{-\sigma R}$ and hence the same type of estimate for $M_v$.

Suppose, more generally, that $M(R)\to0$ is nonincreasing and satisfies $M(R)\le A M(qR)^m+Ae^{-\kappa R}$, where $0<q<1$, $m>1$, and $qm>1$. Set $N=A^{1/(m-1)}M$, so that $N(R)\le N(qR)^m+Be^{-\kappa R}$. Choose $0<\epsilon<1$ with $\epsilon^{qm-1}\le1/2$, and then take $R_0$ sufficiently large that $N\le\epsilon$ on $[R_0,\infty)$. Put $\sigma=-q\log\epsilon/R_0$. Increasing $R_0$ if necessary ensures $\sigma\le\kappa/2$ and $Be^{-(\kappa-\sigma)R_0/q}\le1/2$. On $[R_0,R_0/q]$, we have $N(R)\le e^{-\sigma R}$. If this estimate holds on the preceding geometric interval, then
\[
 N(R)\le e^{-\sigma R}
 \bigl(e^{-\sigma(qm-1)R}+Be^{-(\kappa-\sigma)R}\bigr)
 \le e^{-\sigma R}.
\]
Induction over these intervals proves the preliminary exponential decay. The derivatives are controlled at the same time by \eqref{eq:direct-tail-logbound}.

We next improve the rates. If $v(t)\le Ce^{-r_vt}$ for $t\ge0$, split the convolution into $z<0$, $0\le z\le t$, and $z>t$. For $t\ge0$, the three terms satisfy
\begin{align*}
 \int_{-\infty}^0G(t-z)v(z)^Q\dd z
 &\le Ce^{-bt}\int_{-\infty}^0e^{bz}v(z)^Q\dd z
 \le Ce^{-bt},\\
 \int_0^tG(t-z)v(z)^Q\dd z
 &\le Ce^{-bt}\int_0^t e^{(b-Qr_v)z}\dd z,\\
 \int_t^\infty G(t-z)v(z)^Q\dd z
 &\le Ce^{at}\int_t^\infty e^{-(a+Qr_v)z}\dd z
 \le Ce^{-Qr_vt}.
\end{align*}
The middle integral equals $t$ if $b=Qr_v$ and otherwise is controlled by the endpoint exponentials. Thus
\[
 u(t)\le C(1+t)e^{-\min\{b,Qr_v\}t}.
\]
Every positive rate strictly below $\min\{b,Qr_v\}$ is consequently admissible for $u$; the converse improvement is analogous. Choose target rates $(\bar r_u,\bar r_v)$ arbitrarily close to $(\rho_u,\rho_v)$ such that
\[
 \bar r_u<b,\quad \bar r_u<Q\bar r_v,\qquad
 \bar r_v<a,\quad \bar r_v<P\bar r_u.
\]
To justify this choice, if $bP>a$ and $aQ>b$, then $(\rho_u,\rho_v)=(b,a)$, and both cross inequalities are already strict at this pair; decreasing both components slightly suffices. If $\rho_u=Q\rho_v$, then $\rho_v<P\rho_u$ because $PQ>1$. Decrease $\rho_v$ slightly and decrease $\rho_u$ by a sufficiently larger amount to ensure $\bar r_u<Q\bar r_v$, retaining the other strict inequality. If $\rho_v=P\rho_u$, interchange the components. The two cross equalities cannot hold simultaneously. Finally, a slight further decrease of the target rates ensures $\bar r_u<b$ and $\bar r_v<a$ as well. Choose $\vartheta<1$ sufficiently close to one that $\bar r_u<\vartheta Q\bar r_v$, $\bar r_v<\vartheta P\bar r_u$, and $\vartheta^2PQ>1$. If necessary, decrease the initial preliminary positive rates below both targets, and set successively
\[
 r_u^{j+1}=\min\{\bar r_u,\vartheta Qr_v^j\},\qquad
 r_v^{j+1}=\min\{\bar r_v,\vartheta Pr_u^{j+1}\}.
\]
The second relation reduces to $r_v^{j+1}=\min\{\bar r_v,\vartheta^2PQ\,r_v^j\}$, so this component reaches its target after finitely many steps; the first then does so as well. This proves \eqref{eq:direct-tail-pre-right}.
\end{proof}

\subsection{Three convolution asymptotics}

\begin{lemma}\label{lem:direct-tail-convolution}
Let $K>0$ be bounded and satisfy
\[
 K(x)\sim k e^{-\kappa x}\quad(x\to+\infty),\qquad
 K(x)\le C e^{\kappa_-x}\quad(x\le0),
 \qquad k,\kappa,\kappa_->0.
\]
Let $h>0$ be $C^1$, with $|h|+|h'|\le Ce^{cy}$ for $y\le0$, where $c>0$, and write $w=K*h$. The leading terms of $w$ and $w'$ in the following three cases are as follows.
\begin{enumerate}
\item If $|h(y)|+|h'(y)|\le Ce^{-(\kappa+g)y}$ for $y\ge0$, where $g>0$, then
\[
 w(t)\sim A e^{-\kappa t},\qquad
 w'(t)\sim-\kappa A e^{-\kappa t},\qquad
 A=k\int_\R e^{\kappa y}h(y)\,\dd y>0.
\]
\item If $0<\rho<\kappa$, $e^{\rho y}h(y)\to d>0$, and $e^{\rho y}h'(y)\to-\rho d$, then
\[
 w(t)\sim dM_K(\rho)e^{-\rho t},\qquad
 w'(t)\sim-\rho dM_K(\rho)e^{-\rho t},\qquad
 M_K(\rho)=\int_\R e^{\rho x}K(x)\,\dd x.
\]
\item If $e^{\kappa y}h(y)\to d>0$ and $e^{\kappa y}h'(y)\to-\kappa d$, then
\[
 w(t)\sim kd\,t e^{-\kappa t},\qquad
 w'(t)\sim-\kappa kd\,t e^{-\kappa t}.
\]
\end{enumerate}
All asymptotics are as $t\to+\infty$. The reflected kernel may be used in place of $K$ when needed.
\end{lemma}

\begin{proof}
In each case $h,h'$ are integrable, and differentiation can be transferred to the source: $w'=K*h'$. In the first case, the global bound $e^{\kappa t}K(t-y)\le Ce^{\kappa y}$ permits dominated convergence for both $h$ and $h'$. Integration by parts then gives $\int e^{\kappa y}h'(y)\,\dd y=-\kappa\int e^{\kappa y}h(y)\,\dd y$, with vanishing boundary terms at both ends.

In the second case, $e^{\rho x}K(x)$ is integrable, while $e^{\rho y}h(y)$ and $e^{\rho y}h'(y)$ are bounded on the whole line. Apply dominated convergence to
\[
 e^{\rho t}(K*h)(t)
 =\int_\R e^{\rho x}K(x)
       \bigl[e^{\rho(t-x)}h(t-x)\bigr]\,\dd x
\]
and to the corresponding identity with $h'$.

In the third case, put $A(x)=e^{\kappa x}K(x)$ and $B(y)=e^{\kappa y}h(y)$. Both are bounded on the line and integrable on the negative half-line, and they tend respectively to $k,d$ at positive infinity. Hence
\[
 e^{\kappa t}w(t)
 =\int_\R A(t-y)B(y)\,\dd y
 =\int_0^t A(t-y)B(y)\,\dd y+O(1).
\]
For large fixed $R$, the factors are uniformly close to $k,d$ on $R\le y\le t-R$, while the two remaining end intervals have total length at most $2R$. Letting first $t\to\infty$ and then $R\to\infty$ gives $t^{-1}e^{\kappa t}w(t)\to kd$. Replacing $B$ by $e^{\kappa y}h'(y)$ proves the derivative asymptotic by the same argument, which also applies to this signed function.
\end{proof}

For $P,Q>1$, the function asymptotics in the next proposition are the cited estimates \eqref{eq:direct-published-physical-tails}, together with regularity at the origin. The convolution lemma supplies the derivative refinement; this does not follow by differentiating a pointwise asymptotic equivalence. The same argument also establishes the function asymptotics when $\min\{P,Q\}\le1$.

\begin{proposition}[Asymptotics at both ends for finite-energy solutions]\label{prop:direct-actual-tails}
Suppose \eqref{eq:direct-tail-system} holds. There are strictly positive constants $A_{u,-},A_{v,-},A_{u,+},A_{v,+}$ such that
\begin{align}
 u(t)&\sim A_{u,-}e^{at},&
 v(t)&\sim A_{v,-}e^{bt}
 &&(t\to-\infty),\label{eq:direct-actual-left}\\
 u(t)&\sim A_{u,+}t^{m_u}e^{-\rho_ut},&
 v(t)&\sim A_{v,+}t^{m_v}e^{-\rho_vt}
 &&(t\to+\infty),\label{eq:direct-actual-right}
\end{align}
where
\[
 \rho_u=\min\{b,aQ\},\quad \rho_v=\min\{a,bP\},\qquad
 m_u=\mathbf1_{\{aQ=b\}},\quad m_v=\mathbf1_{\{bP=a\}}.
\]
Moreover,
\begin{equation}\label{eq:direct-actual-loglimits}
 \lim_{t\to-\infty}\frac{u'}u=a,\quad
 \lim_{t\to-\infty}\frac{v'}v=b,\qquad
 \lim_{t\to+\infty}\frac{u'}u=-\rho_u,\quad
 \lim_{t\to+\infty}\frac{v'}v=-\rho_v.
\end{equation}
\end{proposition}

\begin{proof}
At the left end, \eqref{eq:direct-tail-origin} gives \eqref{eq:direct-actual-left}, with $A_{u,-}=U(0)$ and $A_{v,-}=V(0)$. Smoothness and positivity at the origin also give
\[
 \frac{u'(t)}{u(t)}
 =a+\frac{e^t U'(e^t)}{U(e^t)}\longrightarrow a,
 \qquad
 \frac{v'(t)}{v(t)}
 =b+\frac{e^t V'(e^t)}{V(e^t)}\longrightarrow b.
\]

At the right end, if $aQ>b$ and $bP>a$, both sources again decay strictly faster than their respective kernel rates, and case~1 applies directly. If $bP<a$, then $\rho_v=bP$ and $Q\rho_v=bPQ>b$. Case~1 first gives $u(t)\sim A e^{-bt}$ and $u'(t)\sim-bAe^{-bt}$. Thus $u^P$ and its derivative satisfy the hypotheses of case~2, which gives $v(t)\sim A^P M_{G^*}(bP)e^{-bPt}$ and the corresponding derivative asymptotic. If $aQ<b$, interchange the components. The two strictly slower-tail conditions cannot occur simultaneously, since their product would imply $PQ<1$.

If $bP=a$, then $aQ=bPQ>b$, so the pure exponential asymptotics for $u$ and $u'$ follow first. Applying case~3 to $u^P$ then gives $v(t)\sim k_-A^P t e^{-at}$ and its derivative asymptotic. If $aQ=b$, interchange the components. Both resonances cannot occur simultaneously because $PQ>1$. This exhausts all cases, without requiring $P,Q\ge1$.
\end{proof}

\subsection{Square integrability of the weighted kernel}

\begin{proposition}\label{prop:direct-tail-hs}
For the solutions above,
\begin{equation}\label{eq:direct-tail-hs}
 \iint_{\R^2}G(t-z)^2u(t)^{P-1}v(z)^{Q-1}
 \,\dd z\,\dd t<\infty.
\end{equation}
\end{proposition}

\begin{proof}
The coefficients in Proposition~\ref{prop:direct-actual-tails} are strictly positive, so the leading terms provide two-sided comparisons at each end. These remain applicable when $P-1$ or $Q-1$ is negative. Raising the resonance factors to fixed real powers introduces at most polynomial growth. Write
\[
 D=\rho_u(P-1)+\rho_v(Q-1),\qquad
 R_1=2b+\rho_u(P-1),\qquad R_2=2a+\rho_v(Q-1).
\]
In the three regions where both sources decay faster than the kernels, where $bP<a$, and where $aQ<b$, respectively, $D$ equals
\[
 bP-a+aQ-b,\qquad b(PQ-1),\qquad a(PQ-1),
\]
and is strictly positive. The same formulas give strict positivity at the resonance parameters. If $\rho_u=b$, then $R_1=b(P+1)>0$; if $\rho_u=aQ<b$, then
\[
 R_1=a(PQ-1)+a+b+(b-aQ)>0.
\]
Similarly, if $\rho_v=a$, then $R_2=a(Q+1)>0$; otherwise,
\[
 R_2=b(PQ-1)+a+b+(a-bP)>0.
\]

We check integrability over the entire $(t,z)$-plane. Modifying the leading-term models near zero to be positive and bounded extends the two-sided comparisons to each complete half-line. It therefore suffices to divide the plane into the six regions below. In the table, $s,d,x,y\ge0$. Each exponential factor, multiplied by a common constant and at most $(1+s+d)^N$ or $(1+x+y)^N$, bounds the integrand from above.
\begin{center}
\begin{tabular}{lll}
\toprule
Region & Coordinates & Exponential factor\\
\midrule
$0\ge t\ge z$ & $t=-s,\ z=-s-d$ & $e^{-2\alpha s-nd}$\\
$0\ge z\ge t$ & $t=-s-d,\ z=-s$ & $e^{-2\alpha s-nd}$\\
$t\ge z\ge0$ & $t=s+d,\ z=s$ & $e^{-Ds-R_1d}$\\
$z\ge t\ge0$ & $t=s,\ z=s+d$ & $e^{-Ds-R_2d}$\\
$t\ge0\ge z$ & $t=x,\ z=-y$ & $e^{-R_1x-ny}$\\
$z\ge0\ge t$ & $t=-x,\ z=y$ & $e^{-nx-R_2y}$\\
\bottomrule
\end{tabular}
\end{center}
For example, in the first row the three factors satisfy
\[
 G(t-z)^2\le Ce^{-2bd},\quad
 u(t)^{P-1}\le Ce^{-a(P-1)s},\quad
 v(z)^{Q-1}\le Ce^{-b(Q-1)(s+d)}.
\]
Their product is $Ce^{-2\alpha s-nd}$, since $a(P-1)+b(Q-1)=2\alpha$ and $b(Q+1)=n$. In the third row the same multiplication yields
\[
 e^{-2bd-\rho_u(P-1)(s+d)-\rho_v(Q-1)s}
 =e^{-Ds-R_1d}.
\]
The other four rows follow by interchanging the coordinates or substituting $z=-y$ or $t=-x$. All exponential coefficients are strictly positive, so these upper bounds are integrable; a slight decrease in the positive exponents absorbs the polynomial factors. The kernel is bounded near the origin, and $u,v$ have positive lower bounds on every compact interval, so the central region is also integrable. This proves \eqref{eq:direct-tail-hs}.
\end{proof}

\section{The first two eigenvalues of a totally positive kernel}\label{sec:spectral}
We use compact self-adjoint operators whose quadratic forms are
nonnegative. This condition is distinct from pointwise positivity of
their integral kernels. All kernels here are real.
Self-adjointness makes eigenvalues real, and taking real and imaginary
parts shows that each complex eigenspace is the complexification of
its real eigenspace. It therefore suffices to work with real eigenfunctions.

\begin{lemma}[Perron--Jentzsch]\label{lem:perron}
Let \(A\) be a nonzero compact self-adjoint nonnegative integral operator
on \(L^2(X,\mu)\), where the underlying measure space is \(\sigma\)-finite.
If its kernel is strictly positive almost everywhere, then its largest
eigenvalue \(\lambda_1>0\) is simple and has an almost everywhere
strictly positive unit eigenfunction.
\end{lemma}
\begin{proof}
This is the compact self-adjoint case of the Perron--Jentzsch theorem;
see Drnov\v{s}ek \cite[Theorem~2, p.~547]{Drnovsek}.
Indeed, for every nonzero \(f\ge0\), strict positivity of the kernel and
Fubini's theorem give \(Af>0\) almost everywhere. Thus \(A\) is
positivity improving and, in particular, irreducible. Since \(A\ne0\)
is self-adjoint and nonnegative,
\(r(A)=\|A\|>0\); compactness makes this spectral value a pole of the
resolvent. All hypotheses of the cited theorem are satisfied.
\end{proof}

\begin{lemma}[The second eigenvalue and one sign change]\label{lem:second}
Let \(K\) be a compact self-adjoint nonnegative integral operator on
\(L^2(\R)\). Suppose its kernel \(k(t,x)\) is strictly positive almost
everywhere and
\begin{equation}\label{eq:strictdet}
 \det[k(t_i,x_j)]_{i,j=1}^2>0
 \quad\text{for almost every }t_1<t_2,\ x_1<x_2.
\end{equation}
List the eigenvalues, with multiplicity and zeros appended if necessary,
as \(\lambda_1\ge\lambda_2\ge\lambda_3\ge\cdots\ge0\). Then
\begin{equation}\label{eq:spectralgap}
 \lambda_1>\lambda_2>\lambda_3,\qquad \lambda_2>0.
\end{equation}
The first two real eigenfunctions can be chosen so that \(\phi_1>0\)
almost everywhere and, for almost every \(t<x\),
\begin{equation}\label{eq:cheby}
 \phi_1(t)\phi_2(x)-\phi_2(t)\phi_1(x)>0.
\end{equation}
The second eigenfunction changes sign exactly once. Conversely, every
nonzero eigenfunction with exactly one sign change belongs to the
second eigenspace.
\end{lemma}
\begin{proof}
This is the order-two oscillation principle; see
Kushel and Zabreiko \cite[Theorem~6.2]{KushelZabreiko} for a general
operator formulation. In the present self-adjoint setting, it follows
from two applications of Lemma~\ref{lem:perron}.
On \(\Omega=\{(t_1,t_2):t_1<t_2\}\), let \(K^{[2]}\) be the
restriction of \(K\otimes K\) to antisymmetric functions, identified
isometrically with \(L^2(\Omega)\). Its kernel is
\(\det[k(t_i,x_j)]_{i,j=1}^2\), so \(K^{[2]}\) is nonzero,
compact, self-adjoint, nonnegative, and positivity improving.
The spectral theorem gives an orthonormal eigenbasis for \(K\);
the corresponding antisymmetric products diagonalize \(K^{[2]}\)
with eigenvalues \(\lambda_i\lambda_j\), \(i<j\).
Since \(K^{[2]}\ne0\), we have \(\lambda_2>0\).
Lemma~\ref{lem:perron} applied to \(K\) makes \(\lambda_1\)
simple, hence \(\lambda_1>\lambda_2\). Applied to \(K^{[2]}\),
it makes its largest eigenvalue \(\lambda_1\lambda_2\) simple;
therefore \(\lambda_2>\lambda_3\).

Choose the first unit eigenfunction \(\phi_1>0\) almost everywhere.
The function
\(\phi_1(t)\phi_2(x)-\phi_2(t)\phi_1(x)\) is a Perron
eigenfunction of \(K^{[2]}\). Choosing the sign of \(\phi_2\)
makes it strictly positive almost everywhere on \(t<x\), proving
\eqref{eq:cheby}.

Set \(q=\phi_2/\phi_1\), so that \(q(t)<q(x)\) for almost every
\(t<x\). Fubini's theorem implies that the essential supremum of
\(q\) on any interval to the left of another interval is at most
its essential infimum on the latter. Since \(q\) is finite almost
everywhere,
\[
 q_0(t)=\operatorname*{ess\,sup}_{x<t}q(x)
\]
is finite for every finite \(t\) and is nondecreasing.
The same order relation gives \(q_0(t)\le q(t)\le q_0(t+)\)
almost everywhere. A monotone function has at most countably many
jumps, so \(q=q_0\) almost everywhere.
The strict inequality in \eqref{eq:cheby} excludes constancy on any
interval of positive length. Finally,
\(\int\phi_1\phi_2=0\) and \(\phi_1>0\) imply that \(\phi_2\)
has both signs on sets of positive measure. Hence \(q_0\), and therefore
\(\phi_2\), changes sign exactly once, from negative to positive.

Conversely, let \(\xi\) be an eigenfunction with exactly one sign change.
It cannot belong to the first eigenspace. If it does not belong to the
second either, self-adjointness yields
\begin{equation}\label{eq:twoorth}
 \int_\R\xi\phi_1\dd t=\int_\R\xi\phi_2\dd t=0.
\end{equation}
Let \(t_0\) separate its two signs. By monotonicity of \(q\), choose
a finite \(c\) between its essential supremum to the left of \(t_0\)
and its essential infimum to the right. Then
\[
 \xi(t)\phi_1(t)\bigl(q(t)-c\bigr)
\]
has the same sign on both sides of \(t_0\) and is nonzero on a set of
positive measure. It is integrable, since it equals
\(\xi\phi_2-c\xi\phi_1\). Its integral cannot vanish, contradicting
\eqref{eq:twoorth}.
\end{proof}

\begin{remark}\label{rem:order}
Exactly one sign change means that a finite point separates two
opposite fixed signs, with a nonzero part on each side; zero points
and intervals are allowed. Lemma~\ref{lem:second} asserts simplicity
only of the first two eigenvalues. No simplicity of higher eigenvalues
is assumed or used.
\end{remark}

\section{Comparison at large relative scales}\label{sec:translation}
Fix two positive finite-energy solutions \((u_0,v_0)\), \((u_1,v_1)\)
as in Theorem~\ref{thm:allEL}. Their power integrals need not agree.
By Proposition~\ref{prop:direct-actual-tails}, the parameter-dependent quantities
\begin{equation}\label{eq:ratesused}
 \kappa_u=\min\{b,aQ\}>0,\qquad
 m_u=\begin{cases}1,&aQ=b,\\0,&aQ\ne b,\end{cases}
\end{equation}
satisfy, for each \(i\),
\begin{align}
 u_i(t)&=A_{i,-}e^{at}(1+o(1)),&
 \frac{u_i'(t)}{u_i(t)}&\longrightarrow a
 &&(t\to-\infty),\label{eq:leftused}\\
 u_i(t)&=A_{i,+}e^{-\kappa_ut}t^{m_u}(1+o(1)),&
 \frac{u_i'(t)}{u_i(t)}&\longrightarrow-\kappa_u
 &&(t\to+\infty),\label{eq:rightused}
\end{align}
where \(A_{i,\pm}>0\). The second components satisfy analogous formulas.

\begin{lemma}[One sign change at large translations]\label{lem:largeshift}
For all sufficiently large \(s>0\), the function
\begin{equation}\label{eq:ds}
 d_s(t)=u_1(t+s)-u_0(t)
\end{equation}
changes sign exactly once, from positive to negative.
\end{lemma}
\begin{proof}
Write \(\kappa=\kappa_u\), \(m=m_u\).
By \eqref{eq:leftused}--\eqref{eq:rightused}, choose \(R>1\) so that
\begin{equation}\label{eq:oppositetails}
 (\log u_0)'(t)\ge\frac a2\quad(t\le-R),\qquad
 (\log u_1)'(t)\le-\frac\kappa2\quad(t\ge R).
\end{equation}
Let \(s>2R\). If \(t\le-s+R\), then \(t+s\le R\) and \(t\le-R\).
Positivity, continuity, and the left asymptotics give \(c_1,C_0>0\)
such that
\[
 u_1(z)\ge c_1e^{az}\quad(z\le R),\qquad
 u_0(t)\le C_0e^{at}\quad(t\le-R).
\]
Consequently, for all sufficiently large \(s\),
\begin{equation}\label{eq:leftcompare}
 \frac{u_1(t+s)}{u_0(t)}\ge\frac{c_1}{C_0}e^{as}>1
 \quad(t\le-s+R).
\end{equation}
On the other hand, the right asymptotics give \(c_0,C_1>0\) such that
\[
 u_0(t)\ge c_0e^{-\kappa t}(1+t_+)^m\quad(t\ge-R),\qquad
 u_1(z)\le C_1e^{-\kappa z}(1+z_+)^m\quad(z\ge R),
\]
where \(t_+=\max\{t,0\}\). Thus, uniformly for \(t\ge-R\),
\begin{equation}\label{eq:rightcompare}
 \frac{u_1(t+s)}{u_0(t)}
 \le\frac{C_1}{c_0}e^{-\kappa s}
 \left(\frac{1+(t+s)_+}{1+t_+}\right)^m
 \le\frac{C_1}{c_0}e^{-\kappa s}(1+s)^m<1
\end{equation}
when \(s\) is large.
In the remaining interval \(-s+R\le t\le-R\),
\eqref{eq:oppositetails} gives
\[
 \frac{\dd}{\dd t}\log\frac{u_1(t+s)}{u_0(t)}
 \le-\frac{a+\kappa}{2}<0.
\]
Together with \eqref{eq:leftcompare} and \eqref{eq:rightcompare}, this
shows that the ratio crosses one exactly once.
\end{proof}
This argument uses only the two logarithmic derivative limits and
the common exponential and polynomial tail orders; no global
unimodality of the cylindrical profiles is needed.

\section{Secant operators and persistence of the second eigenvalue}\label{sec:secant}
If \(u_1(\cdot+s_0)=u_0\) for some \(s_0\), the second equation in
\eqref{eq:cylEL} also gives \(v_1(\cdot+s_0)=v_0\), and uniqueness
follows. We henceforth assume, for a contradiction, that the profiles
belong to distinct translation orbits. Thus
\begin{equation}\label{eq:distinct}
 d_s\not\equiv0\qquad(s\in\R).
\end{equation}

\subsection{An exact eigenvalue equation for finite differences}
For \(\gamma>0\) and \(x,y>0\), define
\begin{equation}\label{eq:secantdef}
 D_\gamma(x,y)=
 \gamma\int_0^1\bigl((1-\theta)y+\theta x\bigr)^{\gamma-1}\dd\theta.
\end{equation}
This positive secant slope satisfies
\(x^\gamma-y^\gamma=D_\gamma(x,y)(x-y)\).
Set
\begin{align}
 e_s(t)&=v_1(t+s)-v_0(t),\notag\\
 A_s(t)&=D_Q(v_1(t+s),v_0(t)),&
 B_s(t)&=D_P(u_1(t+s),u_0(t)).
 \label{eq:weights}
\end{align}
Subtracting the two systems gives
\begin{equation}\label{eq:diffEL}
 d_s=\Tg(A_se_s),\qquad e_s=\Tg^*(B_sd_s).
\end{equation}
On the fixed Hilbert space \(L^2(\R,\dd t)\), define
\begin{equation}\label{eq:Ls}
 (L_sh)(t)=\sqrt{B_s(t)}
 \int_\R G(t-z)\sqrt{A_s(z)}h(z)\dd z,\qquad K_s=L_sL_s^*.
\end{equation}
Boundedness is established directly for the integral kernel below;
the possibly unbounded multipliers need not be treated separately.

\begin{proposition}[Compactness and the secant eigenfunction]\label{prop:compactsec}
For each finite \(s\), \(L_s\) is Hilbert--Schmidt. Hence \(K_s\) is
compact, self-adjoint, and nonnegative. Moreover,
\begin{equation}\label{eq:fixedmode}
 \xi_s=\sqrt{B_s}\,d_s\in L^2(\R)\setminus\{0\},\qquad K_s\xi_s=\xi_s.
\end{equation}
The kernel of \(K_s\) is strictly positive and strictly totally positive
of order two, almost everywhere.
\end{proposition}
\begin{proof}
The positive leading coefficients in the two tail expansions imply
that, for fixed \(s\), both \(u_1(t+s)/u_0(t)\) and
\(v_1(t+s)/v_0(t)\) are bounded above and bounded away from zero.
Consequently,
\begin{equation}\label{eq:weightcompare}
 B_s(t)\asymp u_0(t)^{P-1},\qquad A_s(t)\asymp v_0(t)^{Q-1}.
\end{equation}
This also holds for negative \(P-1\) or \(Q-1\).
Proposition~\ref{prop:direct-tail-hs} yields
\begin{equation}\label{eq:secHS}
 \|L_s\|_{\HS}^2
 =\iint_{\R^2}G(t-z)^2B_s(t)A_s(z)\dd t\dd z<\infty.
\end{equation}
Furthermore,
\begin{equation}\label{eq:dL2}
 \|\xi_s\|_2^2=
 \int_\R
 \bigl(u_1(t+s)^P-u_0(t)^P\bigr)
 \bigl(u_1(t+s)-u_0(t)\bigr)\dd t<\infty,
\end{equation}
by \(u_i\in L^{P+1}\) and H\"older's inequality.
Similarly, \(\eta_s=\sqrt{A_s}e_s\in L^2\).
Thus \eqref{eq:diffEL} becomes
\(\xi_s=L_s\eta_s\), \(\eta_s=L_s^*\xi_s\), proving
\eqref{eq:fixedmode}; nonvanishing follows from \eqref{eq:distinct}.

The kernel is
\begin{equation}\label{eq:Ks}
 k_s(t,x)=\sqrt{B_s(t)B_s(x)}
 \int_\R G(t-z)A_s(z)G(x-z)\dd z.
\end{equation}
It is finite and strictly positive almost everywhere.
Proposition~\ref{prop:direct-kernel-gram} gives strict total positivity
of order two, also when \(\alpha=2\). The required absolute
integrability follows from \eqref{eq:secHS} and Cauchy--Schwarz.
\end{proof}

\subsection{Continuity in the relative scale}
\begin{proposition}[Norm continuity]\label{prop:continuity}
The map \(s\mapsto L_s\) is continuous in Hilbert--Schmidt norm,
and \(s\mapsto K_s\) is continuous in operator norm.
\end{proposition}
\begin{proof}
By \eqref{eq:direct-tail-logbound}, there are finite \(C_u,C_v\) with
\[
 |u_1'/u_1|\le C_u,\qquad |v_1'/v_1|\le C_v.
\]
Hence, for all \(t,s,h\in\R\),
\begin{equation}\label{eq:logratio}
 e^{-C_u|h|}\le\frac{u_1(t+s+h)}{u_1(t+s)}\le e^{C_u|h|},
\end{equation}
and similarly for \(v_1\).
The same bounds hold for the ratio of the corresponding convex
combinations with the fixed positive number \(u_0(t)\).
Raising them to the power \(P-1\) and integrating in
\eqref{eq:secantdef} gives
\begin{align}
 e^{-|P-1|C_u|h|}
 &\le B_{s+h}(t)/B_s(t)\le e^{|P-1|C_u|h|},\notag\\
 e^{-|Q-1|C_v|h|}
 &\le A_{s+h}(t)/A_s(t)\le e^{|Q-1|C_v|h|}.
 \label{eq:seccontratio}
\end{align}
Both \(\beta_{s,h}=(B_{s+h}/B_s)^{1/2}\) and
\(\gamma_{s,h}=(A_{s+h}/A_s)^{1/2}\) tend to one in \(L^\infty\).
Writing \(M_w\) for multiplication by \(w\), we have
\[
 L_{s+h}=M_{\beta_{s,h}}L_sM_{\gamma_{s,h}}.
\]
The ideal property of the Hilbert--Schmidt norm gives
\(\|L_{s+h}-L_s\|_{\HS}\to0\).
Finally,
\[
 \|L_{s+h}L_{s+h}^*-L_sL_s^*\|_{\op}
 \le(\|L_{s+h}\|_{\op}+\|L_s\|_{\op})\|L_{s+h}-L_s\|_{\op},
\]
which proves the second assertion.
\end{proof}

\subsection{Persistence of the second spectral position}
List the eigenvalues of \(K_s\), with multiplicity, as
\(\Lambda_1(s)\ge\Lambda_2(s)\ge\Lambda_3(s)\ge\cdots\ge0\).
Proposition~\ref{prop:compactsec} and Lemma~\ref{lem:second} give
\begin{equation}\label{eq:allgaps}
 \Lambda_1(s)>\Lambda_2(s)>\Lambda_3(s),\qquad
 \Lambda_2(s)>0\quad(s\in\R).
\end{equation}
Equation~\eqref{eq:fixedmode} makes \(1\) an eigenvalue for every \(s\).
Its spectral position is first identified at large positive translations
and then propagated by continuity.

\begin{proposition}[The secant difference is a second eigenfunction]\label{prop:secondall}
Under \eqref{eq:distinct}, for every \(s\in\R\),
\begin{equation}\label{eq:secondall}
 \Lambda_2(s)=1>\Lambda_3(s).
\end{equation}
Consequently, \(d_s\) changes sign exactly once at every relative scale.
\end{proposition}
\begin{proof}
For sufficiently large \(s\), Lemma~\ref{lem:largeshift} gives exactly
one sign change of \(d_s\), and hence of \(\xi_s=\sqrt{B_s}d_s\).
Equations~\eqref{eq:fixedmode} and Lemma~\ref{lem:second} imply
\(\Lambda_2(s)=1\) there. Set
\[
 E=\{s\in\R:\Lambda_2(s)=1\}.
\]
The min--max principle, applied to finite-dimensional compressions
\cite[Chapter~I, \S6.10, pp.~60--61]{Kato}, gives
\[
 |\Lambda_j(s)-\Lambda_j(s_0)|\le\|K_s-K_{s_0}\|_{\op}.
\]
Thus the eigenvalues are continuous and \(E\) is nonempty and closed.
If \(s_0\in E\), then \eqref{eq:allgaps} gives
\(\Lambda_1(s_0)>1>\Lambda_3(s_0)\).
These inequalities persist near \(s_0\).
The eigenvalue \(1\), which exists at every nearby \(s\), can then
occupy only the second spectral position. Hence \(E\) is open as well.
Connectedness of \(\R\) gives \(E=\R\), and a final application of
Lemma~\ref{lem:second} proves the sign-change assertion.
\end{proof}

\section{Classification of all positive finite-energy solutions}\label{sec:allel-proof}

No equality of the power integrals of the two solutions is assumed in this section. The preceding section established that, unless two solutions are common translates of one another, their nonzero difference has exactly one sign change at every relative translation. The orientation of this sign change completes the classification.

\begin{proof}[Proof of Theorem~\ref{thm:allEL}]
Let the two positive solutions be $(u_0,v_0)$ and $(u_1,v_1)$, and set
\[
 d_s(t)=u_1(t+s)-u_0(t).
\]
If $d_{s_*}\equiv0$ for some $s_*\in\R$, the second cylindrical equation gives
\[
 v_1(t+s_*)=\Tg^*\bigl(u_1(\cdot+s_*)^P\bigr)(t)
           =\Tg^*(u_0^P)(t)=v_0(t),
\]
which proves the theorem. Suppose instead that no such $s_*$ exists. Then $d_s\not\equiv0$ for every $s\in\R$, and Proposition~\ref{prop:secondall} implies that each $d_s$ has exactly one sign change. This spectral conclusion was initially obtained almost everywhere. Since $d_s$ is continuous and the weight $\sqrt{B_s}$ converting it into an eigenfunction is everywhere positive, every strictly positive or strictly negative point has a neighborhood of the same sign. Consequently, the almost-everywhere sign-change statement also constrains the ordering of all points where the sign is strict.

We first determine the orientation for large translations of either sign. Lemma~\ref{lem:largeshift} shows that, for sufficiently large positive $s$, $d_s$ is positive to the left of its sign change and negative to the right. For sufficiently large $h>0$, put $s=-h$ and $z=t-h$. Then
\[
 d_{-h}(t)=u_1(z)-u_0(z+h).
\]
Interchanging the solutions and applying Lemma~\ref{lem:largeshift} again shows that $u_0(z+h)-u_1(z)$ is positive on the left and negative on the right. Negation reverses these signs, and $t\mapsto z=t-h$ preserves order. Thus, for sufficiently large negative $s$, $d_s$ is negative on the left and positive on the right.

To prove that the orientation cannot change at an intermediate parameter, define
\begin{align}
 \mathcal O_+
 &=\{s\in\R:\text{there exist }t_1<t_2,\ 
                  d_s(t_1)>0>d_s(t_2)\},\label{eq:orientationplus}\\
 \mathcal O_-
 &=\{s\in\R:\text{there exist }t_1<t_2,\ 
                  d_s(t_1)<0<d_s(t_2)\}.\label{eq:orientationminus}
\end{align}
Because every $d_s$ has exactly one sign change, these sets are disjoint and cover $\R$. Both are open: at any member, the two fixed observation points in the definition retain their strict signs under a small change of $s$. The large-translation comparisons make both sets nonempty, contradicting connectedness of $\R$. Hence $d_{s_*}\equiv0$ for some $s_*$, and the other equation identifies the second components.

Only fixed finite observation points enter the openness argument; no continuous choice or transversality of the sign-change point is required. Zeros and zero intervals, or a sign-change point tending to infinity, cause no difficulty. No equality of masses between the two solutions has been used.

For the physical-space solutions in the final assertion of the theorem, Hang's theorem, cited above, makes both components radial about a common center. Translate the center of each solution to the origin and apply \eqref{eq:EF}. The finite-energy conditions become precisely the cylindrical integrability conditions used above. Under the inverse transformation, the common cylindrical translation is the common dilation
\((U,V)(x)\mapsto(\rho^aU(\rho x),\rho^bV(\rho x))\)
for a single $\rho>0$. Restoring the common centers gives the physical-space classification.
\end{proof}

\begin{proof}[Proof of Theorem~\ref{thm:main}]
Apply \eqref{eq:amplitude} to $(f_*,g_*)$ and to any normalized extremizing pair $(f,g)$. Both resulting potential pairs solve \eqref{eq:physicalEL} and have finite energy, so Theorem~\ref{thm:allEL} identifies them up to a common translation and dilation. Recovering the sources gives dilation powers $aP=n/p$ and $bQ=n/r$; the amplitude constants agree because both pairs have multiplier $C_\sharp$. This proves \eqref{eq:uniqueness}. Conversely, that transformation preserves both norms and the functional, since $n/p+n/r-2n+\lambda=0$.
\end{proof}

\begin{corollary}[Normalized positive radial stationary pairs are extremizing]\label{cor:stationarysharp}
Let $f,g\ge0$ be radial about the same center, with $\|f\|_p=\|g\|_r=1$, and suppose that, for some $\ell>0$,
\begin{equation}\label{eq:stationaryell}
 Tg=\ell f^{p-1},\qquad Tf=\ell g^{r-1}.
\end{equation}
Then $\ell=C_\sharp$, and $(f,g)$ is a normalized extremizing pair.
\end{corollary}
\begin{proof}
Strict positivity of the Riesz kernel implies $f,g>0$. Multiplying the first equation by $f$ and integrating gives $\mathcal J(f,g)=\ell$; the second equation gives the same value.

Using the amplitude normalization above, set
\[
 U=\ell^{-(Q+1)/(PQ-1)}f^{p-1},\qquad
 V=\ell^{-(P+1)/(PQ-1)}g^{r-1}.
\]
Direct substitution into \eqref{eq:stationaryell} gives $U=T(V^Q)$ and $V=T(U^P)$, with
\begin{equation}\label{eq:ellmass}
 \int_{\R^n}U^{P+1}\dd x
 =\int_{\R^n}V^{Q+1}\dd x
 =\ell^{-(P+1)(Q+1)/(PQ-1)}
 =\ell^{-n/\alpha}.
\end{equation}
Apply the same normalization to an existing extremizing pair, replacing $\ell$ by $C_\sharp$. Theorem~\ref{thm:allEL} relates the two potential pairs by translations and dilations, which preserve their power integrals because $a(P+1)=b(Q+1)=n$. Thus \eqref{eq:ellmass} gives $\ell^{-n/\alpha}=C_\sharp^{-n/\alpha}$, and hence $\ell=C_\sharp$.
\end{proof}

\section{Stability of an individual solution and the radial linearized equations}\label{sec:single-spectrum}

Fix any positive finite-energy cylindrical solution
\[
 u=\Tg(v^Q),\qquad v=\Tg^*(u^P).
\]
The spectral conclusions below do not assume in advance that this solution is extremizing. We first construct probability operators from the integral equations themselves and then compare the solution with its nonzero translates.

\subsection{Conditional probabilities determined by the equations}

Multiplying each equation by the corresponding source and applying Tonelli's theorem gives
\begin{equation}\label{eq:singlemass}
 S:=\int_\R u^{P+1}\dd t
   =\iint_{\R^2}u(t)^PG(t-z)v(z)^Q\dd t\dd z
   =\int_\R v^{Q+1}\dd z.
\end{equation}
All these integrals are finite and positive. Define three probability measures by
\begin{align}
 \dd\mu(t)&=S^{-1}u(t)^{P+1}\dd t,&
 \dd\nu(z)&=S^{-1}v(z)^{Q+1}\dd z,\label{eq:singlemarginals}\\
 \dd\pi(t,z)&=S^{-1}u(t)^PG(t-z)v(z)^Q\dd t\dd z.
 \label{eq:singlejoint}
\end{align}
The first equation identifies the $t$-marginal of $\pi$ as $\mu$, and the second identifies its $z$-marginal as $\nu$.

The conditional probabilities at fixed $t$ and fixed $z$ are, respectively,
\begin{equation}\label{eq:conditionallaws}
 \dd R_t(z)=\frac{G(t-z)v(z)^Q}{u(t)}\dd z,\qquad
 \dd R_z^*(t)=\frac{G(t-z)u(t)^P}{v(z)}\dd t.
\end{equation}
Both measures have total mass one. This uses the two equations, beyond positivity of the kernel alone. For real-valued functions, define
\begin{equation}\label{eq:conditionaloperators}
 (\mathcal C k)(t)=\int_\R k(z)\dd R_t(z)
                 =\frac{\Tg(v^Qk)(t)}{u(t)},\qquad
 (\mathcal C^*h)(z)=\int_\R h(t)\dd R_z^*(t)
                 =\frac{\Tg^*(u^Ph)(z)}{v(z)}.
\end{equation}
As shown next, $\mathcal C^*$ is indeed the adjoint between the probability $L^2$ spaces.

\begin{proposition}[Contraction of conditional averages and conjugacy of the integral operators]\label{prop:conditionalmarkov}
For every $1\le q\le\infty$,
\begin{equation}\label{eq:conditionalcontraction}
 \|\mathcal Ck\|_{L^q(\mu)}\le\|k\|_{L^q(\nu)},\qquad
 \|\mathcal C^*h\|_{L^q(\nu)}\le\|h\|_{L^q(\mu)}.
\end{equation}
Both operators preserve nonnegativity and constants, and satisfy
\[
 \int\mathcal Ck\dd\mu=\int k\dd\nu,\qquad
 \int\mathcal C^*h\dd\nu=\int h\dd\mu.
\]
For $q=2$, they are adjoints.

Define the isometric isomorphisms
\[
 \mathscr U h=S^{-1/2}u^{(P+1)/2}h,\qquad
 \mathscr V k=S^{-1/2}v^{(Q+1)/2}k
\]
from $L^2(\mu)$ and $L^2(\nu)$, respectively, onto $L^2(\R,\dd t)$. Then
\begin{equation}\label{eq:singleconjugacies}
 \mathcal L=\mathscr U\mathcal C\mathscr V^{-1},\qquad
 \mathcal L^*=\mathscr V\mathcal C^*\mathscr U^{-1},
\end{equation}
where
\begin{equation}\label{eq:singleL}
 (\mathcal L\chi)(t)
 =u(t)^{(P-1)/2}\int_\R
       G(t-z)v(z)^{(Q-1)/2}\chi(z)\dd z.
\end{equation}
The operator $\mathcal L$ is Hilbert--Schmidt and $\|\mathcal L\|_{\op}=1$.
\end{proposition}
\begin{proof}
On the joint probability space $(\R^2,\pi)$, the operators $\mathcal C$ and $\mathcal C^*$ are conditional expectations onto the two coordinate sigma-fields. Their $L^q$ contraction property is the standard consequence of conditional Jensen's inequality; see Kallenberg \cite[Chapter~6, p.~109]{Kallenberg}. This proves \eqref{eq:conditionalcontraction} for $1\le q<\infty$; the $L^\infty$ case follows from positivity and preservation of constants. The marginal identities above identify the coordinate norms with $L^q(\mu)$ and $L^q(\nu)$ and also give the stated integral identities.

For $h\in L^2(\mu)$ and $k\in L^2(\nu)$,
\[
 \iint|h(t)k(z)|\dd\pi(t,z)
 \le\left(\int|h|^2\dd\mu\right)^{1/2}
     \left(\int|k|^2\dd\nu\right)^{1/2}<\infty.
\]
Fubini's theorem therefore gives
\begin{equation}\label{eq:conditionaladjoint}
 \langle h,\mathcal Ck\rangle_{L^2(\mu)}
 =\iint h(t)k(z)\dd\pi(t,z)
 =\langle\mathcal C^*h,k\rangle_{L^2(\nu)}.
\end{equation}

The isometric isomorphism properties follow immediately from \eqref{eq:singlemarginals}. Substituting $\mathscr V^{-1}\chi=S^{1/2}v^{-(Q+1)/2}\chi$ gives
\[
 \mathscr U\mathcal C\mathscr V^{-1}\chi
 =u^{(P-1)/2}\Tg\bigl(v^{(Q-1)/2}\chi\bigr),
\]
which is \eqref{eq:singleL}; the other identity follows by the same computation or by taking adjoints. Thus \eqref{eq:singleL} defines a bounded integral operator even when $P<1$ or $Q<1$, without requiring its individual multiplication factors to be bounded operators.

Proposition~\ref{prop:direct-tail-hs} gives the Hilbert--Schmidt property of $\mathcal L$. Since $\mathcal C$ is an $L^2$ contraction, $\|\mathcal L\|_{\op}\le1$. On the other hand, set
\[
 \phi=S^{-1/2}u^{(P+1)/2}=\mathscr U\one,\qquad
 \psi=S^{-1/2}v^{(Q+1)/2}=\mathscr V\one.
\]
Both are unit vectors, and preservation of constants gives
\begin{equation}\label{eq:singleperronvectors}
 \mathcal L\psi=\phi,\qquad \mathcal L^*\phi=\psi.
\end{equation}
The operator norm is therefore exactly one.
\end{proof}

\begin{remark}
The Markov property here refers to preservation of nonnegativity and constants and to the averaging structure on the probability spaces above. For $1<q<\infty$, strict convexity also shows that $\|\mathcal Ck\|_q=\|k\|_q$ forces $k$ to be constant almost everywhere. Indeed, equality after integration requires equality in conditional Jensen almost everywhere, and each $R_t$ is equivalent to $\nu$, since both have everywhere strictly positive densities. This equality characterization alone does not give a uniform norm bound strictly below one on the centered subspace; the spectral argument below determines the required value.
\end{remark}

\subsection{Determination of the second eigenvalue by self-translation}

\begin{theorem}[Radial spectrum of an individual positive solution]\label{thm:radialspectrum}
Let
\[
 \mathcal B=\mathcal L\mathcal L^*
\]
act on $L^2(\R,\dd t)$, and list its eigenvalues with multiplicity as $\lambda_1\ge\lambda_2\ge\lambda_3\ge\cdots\ge0$, appending zeros if necessary. Then
\begin{equation}\label{eq:exactradialspectrum}
 1=\lambda_1>\lambda_2=\frac1{PQ}>\lambda_3.
\end{equation}
The first two eigenspaces are one-dimensional, and
\begin{equation}\label{eq:radialthresholdvector}
 \ker\bigl(\mathcal B-(PQ)^{-1}I\bigr)
 =\operatorname{span}\{u^{(P-1)/2}u'\}.
\end{equation}
\end{theorem}
\begin{proof}
By \eqref{eq:singleperronvectors}, $\mathcal B\phi=\phi$, while $\|\mathcal L\|_{\op}=1$. Thus $\lambda_1(\mathcal B)=1$. The integral kernel is strictly positive, so Lemma~\ref{lem:perron} makes the first eigenspace simple.

Compare $(u,v)$ with $(u(\cdot+s),v(\cdot+s))$, taking only $s>0$. If $u(\cdot+s)=u$ for some nonzero $s$, then $u$ is periodic. A continuous strictly positive periodic function has a positive $L^{P+1}$ integral over one period and hence infinite integral on the whole line, contrary to finite energy. Thus the translated differences are nonzero for every $s>0$.

Let $A_s,B_s,L_s,K_s$ be the secant weights and operators previously defined for these two solutions. Propositions~\ref{prop:compactsec} and~\ref{prop:continuity} apply for every $s>0$, and Lemma~\ref{lem:largeshift} gives exactly one sign change for sufficiently large $s$. To exclude the identically zero difference at $s=0$, repeat the proof of Proposition~\ref{prop:secondall} on the connected interval $(0,\infty)$. The set
\[
 E_+=\{s>0:\Lambda_2(K_s)=1\}
\]
is nonempty and closed. At each of its points, strict order-two total positivity gives $\Lambda_1>1>\Lambda_3$. These strict inequalities persist in a neighborhood, where a nonzero eigenvector with eigenvalue one remains present. Hence $E_+$ is also open, and
\begin{equation}\label{eq:selfshiftsecond}
 \Lambda_2(K_s)=1\qquad(s>0).
\end{equation}

Although the difference vanishes identically at $s=0$, the secant weights and operators have well-defined limits. By the secant definition,
\[
 A_0=Qv^{Q-1},\qquad B_0=Pu^{P-1}.
\]
The relative multiplier estimates in Proposition~\ref{prop:continuity} also apply at zero, giving
\begin{equation}\label{eq:selfshiftlimit}
 L_s\longrightarrow\sqrt{PQ}\,\mathcal L\quad\text{in Hilbert--Schmidt norm},\qquad
 K_s\longrightarrow PQ\,\mathcal B\quad\text{in operator norm}.
\end{equation}
The eigenvalue continuity used in Proposition~\ref{prop:secondall}, together with \eqref{eq:selfshiftsecond}, implies
\[
 PQ\,\lambda_2(\mathcal B)=1.
\]
Moreover, $\mathcal B$ itself has a strictly positive kernel with strict order-two total positivity: apply Proposition~\ref{prop:direct-kernel-gram} with the positive weight $v^{Q-1}$. Lemma~\ref{lem:second} therefore gives simplicity of the second eigenvalue and $\lambda_3<1/(PQ)$. Strictness follows from the limiting operator's own properties, rather than from taking a limit of strict inequalities.

It remains to identify the second eigenfunction. Proposition~\ref{prop:direct-actual-tails} and \eqref{eq:direct-tail-logbound} give
\[
 \zeta=u^{(P-1)/2}u'\in L^2,\qquad
 \|\zeta\|_2^2
 \le\|u'/u\|_\infty^2\int u^{P+1}\dd t<\infty.
\]
This function is nonzero because $u$ is nonzero and tends to zero at both ends. Exponential decay and bounded logarithmic derivatives also give $u^P,v^Q\in W^{1,1}(\R)$; for instance, $(u^P)'=Pu^P(u'/u)$, including when $P<1$. Thus derivatives may be transferred to the sources in the convolutions, yielding
\begin{equation}\label{eq:actualderivativeEL}
 u'=Q\Tg(v^{Q-1}v'),\qquad
 v'=P\Tg^*(u^{P-1}u').
\end{equation}
Setting $\eta=v^{(Q-1)/2}v'\in L^2$, these become
\[
 \zeta=Q\mathcal L\eta,\qquad
 \eta=P\mathcal L^*\zeta.
\]
Hence $\mathcal B\zeta=(PQ)^{-1}\zeta$. The second eigenspace is one-dimensional, proving \eqref{eq:radialthresholdvector}.
\end{proof}

\subsection{Second variation and the function spaces for linearization}

Here the potential variables are $u,v$, and the corresponding source variables are $u^P,v^Q$. Write $h,k$ for the relative source variations, so that the first variations of the sources are $u^Ph,v^Qk$. The first-order source norm constraints are precisely
\[
 \int h\dd\mu=0,\qquad \int k\dd\nu=0.
\]
Denote the centered spaces by
\[
 L^2_0(\mu)=\left\{h\in L^2(\mu):\int h\dd\mu=0\right\},\qquad
 L^2_0(\nu)=\left\{k\in L^2(\nu):\int k\dd\nu=0\right\}.
\]
Since conditional averaging preserves integrals, the operators map the corresponding centered spaces into one another.

By \eqref{eq:singleconjugacies}, $\mathcal C\mathcal C^*$ is isometrically conjugate to $\mathcal B$, with constants corresponding to the first eigenvector $\phi$. Theorem~\ref{thm:radialspectrum} therefore gives the exact centered norm bound
\begin{equation}\label{eq:centeredconditionalnorm}
 \|\mathcal Ck\|_{L^2(\mu)}
 \le(PQ)^{-1/2}\|k\|_{L^2(\nu)},\qquad k\in L^2_0(\nu).
\end{equation}
Consequently, for centered $h,k$,
\begin{align}
 \mathscr Q(h,k)
 &:=\frac1P\|h\|_{L^2(\mu)}^2+
     \frac1Q\|k\|_{L^2(\nu)}^2
       -2\langle h,\mathcal Ck\rangle_{L^2(\mu)}
 \notag\\
 &\ge
 \left(\frac{\|h\|_{L^2(\mu)}}{\sqrt P}
       -\frac{\|k\|_{L^2(\nu)}}{\sqrt Q}\right)^2
 \ge0.\label{eq:sourcehessianpositive}
\end{align}

This is the stability quadratic form under the HLS norm constraints. To verify this directly, continue to use the unnormalized cylindrical kernel $G$, and put
\[
 F=\frac{u^P}{S^{1/p}},\qquad H=\frac{v^Q}{S^{1/r}}.
\]
Then $\|F\|_p=\|H\|_r=1$ and
\[
 J_G(F,H):=\int F\,\Tg H\dd t=S^{-\alpha/n}=:j_0.
\]
For bounded centered $h,k$ and sufficiently small $|\varepsilon|$, define
\[
 F_\varepsilon=\frac{F(1+\varepsilon h)}
                       {\|F(1+\varepsilon h)\|_p},\qquad
 H_\varepsilon=\frac{H(1+\varepsilon k)}
                       {\|H(1+\varepsilon k)\|_r}.
\]
Both sources remain nonnegative. Taylor expansion, $p-1=1/P$, $r-1=1/Q$, and \eqref{eq:conditionaladjoint} give
\begin{equation}\label{eq:sourcevariation}
 J_G(F_\varepsilon,H_\varepsilon)
 =j_0\left(1-\frac{\varepsilon^2}{2}\mathscr Q(h,k)
                +o(\varepsilon^2)\right).
\end{equation}
For example, the first denominator expands as $1+\frac{p-1}{2}\varepsilon^2\int h^2\dd\mu+o(\varepsilon^2)$, while the quadratic cross term of the numerator is $j_0\varepsilon^2\langle h,\mathcal Ck\rangle_\mu$. These terms give \eqref{eq:sourcevariation} directly. Bounded centered functions are dense in the centered $L^2$ spaces, so nonnegativity of $\mathscr Q$ extends to the full Hilbert space stated above. No twice differentiability on an entire $L^p$ neighborhood is inferred from this quadratic-form calculation.

\begin{corollary}[Radial linearized solutions and a strict estimate in orthogonal directions]\label{cor:radialjacobi}
Suppose the potential variation $(\dot u,\dot v)$ belongs to
\begin{equation}\label{eq:jacobidomain}
 \mathcal Y_u\times\mathcal Y_v
 =L^2\bigl(\R,u^{P-1}\dd t\bigr)
       \times L^2\bigl(\R,v^{Q-1}\dd t\bigr)
\end{equation}
and satisfies the linearized equations
\begin{equation}\label{eq:jacobiequations}
 \dot u=Q\Tg(v^{Q-1}\dot v),\qquad
 \dot v=P\Tg^*(u^{P-1}\dot u),
\end{equation}
interpreted through the bounded weighted operators specified below. Then, for some $c\in\R$,
\begin{equation}\label{eq:jacobiscaleonly}
 (\dot u,\dot v)=c(u',v').
\end{equation}

On the centered source space $L^2_0(\mu)\times L^2_0(\nu)$, the nullspace of $\mathscr Q$ is
\[
 \operatorname{span}\{(Pu'/u,Qv'/v)\}.
\]
Equip this product space with the inner product
\[
 \langle(h,k),(\widetilde h,\widetilde k)\rangle_w
 =P^{-1}\langle h,\widetilde h\rangle_\mu+
  Q^{-1}\langle k,\widetilde k\rangle_\nu.
\]
If $(h,k)$ is orthogonal to this one-dimensional nullspace in that inner product, then
\begin{equation}\label{eq:transversequadraticgap}
 \mathscr Q(h,k)\ge
 \bigl(1-\sqrt{PQ\,\lambda_3(\mathcal B)}\bigr)
 \left(P^{-1}\|h\|_\mu^2+Q^{-1}\|k\|_\nu^2\right).
\end{equation}
The coefficient in parentheses is strictly positive.
\end{corollary}
\begin{proof}
Set
\[
 \chi=u^{(P-1)/2}\dot u,\qquad
 \vartheta=v^{(Q-1)/2}\dot v.
\]
The natural domain \eqref{eq:jacobidomain} is exactly the condition $\chi,\vartheta\in L^2$. By \eqref{eq:singleL}, the linearized equations mean precisely
\[
 \chi=Q\mathcal L\vartheta,\qquad
 \vartheta=P\mathcal L^*\chi.
\]
These are identities for bounded operators on the full $L^2$ spaces. Eliminating the second component gives $\mathcal B\chi=(PQ)^{-1}\chi$. By \eqref{eq:radialthresholdvector}, $\dot u=cu'$; substitution into the second equation and \eqref{eq:actualderivativeEL} then give $\dot v=cv'$. Conversely, the theorem's proof already verifies that $(u',v')$ belongs to the natural domain and solves these equations.

The relative source variations and potential variations are related by
\begin{equation}\label{eq:sourcepotentialdirections}
 h=P\dot u/u,\qquad k=Q\dot v/v.
\end{equation}
Thus the potential equations are equivalent to $h=P\mathcal Ck$ and $k=Q\mathcal C^*h$. These relations automatically imply centering: if $m_h=\int h\dd\mu$ and $m_k=\int k\dd\nu$, then $m_h=Pm_k$ and $m_k=Qm_h$. Since $PQ>1$, both vanish.

The quadratic form corresponds, in the inner product $\langle\cdot,\cdot\rangle_w$, to the bounded self-adjoint operator
\[
 \mathscr A(h,k)=(h-P\mathcal Ck,\ k-Q\mathcal C^*h).
\]
By \eqref{eq:sourcehessianpositive}, this operator is nonnegative. The quadratic form of a nonnegative self-adjoint operator vanishes exactly on its kernel. The classification of linearized solutions just proved therefore also identifies the nullspace of $\mathscr Q$.

Finally, the isometric change of variables $(h,k)\mapsto(h/\sqrt P,k/\sqrt Q)$ transforms $\mathscr A$ into
\[
 \begin{pmatrix}
 I&-\sqrt{PQ}\,\mathcal C\\
 -\sqrt{PQ}\,\mathcal C^*&I
 \end{pmatrix}
\]
on the centered spaces. The nonzero singular values of the conditional operator are the square roots of the eigenvalues of $\mathcal B$ with its first eigenvalue removed. The largest is therefore $1/\sqrt{PQ}$ and is simple; all the others are at most $\sqrt{\lambda_3(\mathcal B)}$. For a singular-vector pair, the corresponding two-dimensional matrix has eigenvalues $1\pm\sqrt{PQ}\,\sigma$; the kernels of the conditional operators contribute eigenvalue one. The unique zero eigenvalue arises from the largest singular value, and its paired eigenvalue is two. After removing this one-dimensional nullspace, the entire spectrum is bounded below by $1-\sqrt{PQ\,\lambda_3(\mathcal B)}>0$. This is \eqref{eq:transversequadraticgap}.
\end{proof}

The derivatives $u',v'$ are the first variations under cylindrical translation; in physical variables they are the infinitesimal direction of the common dilation. The preceding arguments establish both the classification of all positive solutions and the strict spectral estimate transverse to this direction. The latter follows from comparison with actual translates, independently of any inference from uniqueness to nondegeneracy.

\section{Normalized alternating iteration and convergence}\label{sec:iteration}

Corollary \ref{cor:stationarysharp} shows that every normalized positive finite-energy Euler--Lagrange solution is an extremizing pair. We now consider alternating iteration from arbitrary nonnegative radially nonincreasing initial data. Each update increases the objective, but convergence of the objective values alone does not imply convergence of the functions. We first identify all possible limiting profiles, then exclude two or more mutually separating translates, and finally prove convergence of the remaining profile at a fixed scale.

\subsection{Kernel normalization and unit norms on the cylinder}

Fix $n\ge3$, $2\le\alpha<n$, and $\lambda=n-\alpha$, and set
\begin{equation}\label{eq:it-basicparameters}
\frac1p+\frac1r=1+\frac{\alpha}{n},\qquad
P=\frac1{p-1},\quad Q=\frac1{r-1},\quad
a=\frac n{P+1},\quad b=\frac n{Q+1}.
\end{equation}
Then $a+b=\lambda$, $aP=b+\alpha$, $bQ=a+\alpha$, and $PQ>1$.
Write $\omega=|\mathbb S^{n-1}|$, and let $T$ denote the unnormalized Riesz integral operator used above. In this part we use the linear isometry
\begin{equation}\label{eq:it-isometry-definition}
(\mathcal U_q h)(t)=\omega^{1/q}e^{nt/q}h(e^t),\qquad
\|\mathcal U_qh\|_{L^q(\R)}=\|h\|_{L^q(\R^n)}.
\end{equation}
Here $h$ is radial and $h(e^t)$ denotes its radial profile.
The associated convolution kernel is
\begin{equation}\label{eq:it-normalized-kernel}
\mathcal K(t)=\omega^{-\alpha/n}G(t),\qquad
Kh=\mathcal K*h,\qquad K^*h=\mathcal K^**h,\quad
\mathcal K^*(t)=\mathcal K(-t).
\end{equation}
Thus, if $F_0=\mathcal U_p f_0$ and $H_0=\mathcal U_r g_0$, then
\begin{equation}\label{eq:it-bilinear-isometry}
\mathcal J(f_0,g_0)
=\iint_{\R^n\times\R^n}|x-y|^{-\lambda}f_0(x)g_0(y)\,\dd x\,\dd y
=\int_\R F_0(t)(KH_0)(t)\,\dd t.
\end{equation}
The factor $\omega^{-\alpha/n}$ equals $\omega^{1-1/p-1/r}$ and is necessary when both norms are fixed at one. Correspondingly,
\[
\|Kh\|_{p'}\le C_\sharp\|h\|_r,\qquad
\|K^*h\|_{r'}\le C_\sharp\|h\|_p.
\]
These inequalities also hold for signed functions, since positivity of the kernel gives $|Kh|\le K|h|$. The positivity, exponential tails, and derivative bounds of Lemma \ref{lem:direct-kernel-factor} remain valid after multiplication by this fixed positive constant.

Fix a unit-norm extremizing pair $(f_*,g_*)$ with common center at the origin, and write
\[
F=\mathcal U_pf_*,\qquad H=\mathcal U_rg_*.
\]
Then
\begin{equation}\label{eq:it-optimal-pair}
\|F\|_p=\|H\|_r=1,\qquad
K^*F=C_\sharp H^{r-1},\qquad KH=C_\sharp F^{p-1}.
\end{equation}
To relate this normalization to the preceding stationary results, suppose that positive unit-norm functions $X,Y$ satisfy $KY=jX^{p-1}$ and $K^*X=jY^{r-1}$. Applying the inverse transforms $x=\mathcal U_p^{-1}X$, $y=\mathcal U_r^{-1}Y$ gives exactly $Ty=jx^{p-1}$ and $Tx=jy^{r-1}$. Corollary \ref{cor:stationarysharp} therefore yields $j=C_\sharp$. Thus the stationary classification used below has precisely the same multiplier as its Euclidean counterpart.

Define auxiliary profiles with coefficient one for the normalized kernel by
\begin{equation}\label{eq:it-auxiliary-profiles}
\mathsf u=C_\sharp^{-(Q+1)/(PQ-1)}F^{1/P},\qquad
\mathsf v=C_\sharp^{-(P+1)/(PQ-1)}H^{1/Q}.
\end{equation}
They satisfy $\mathsf u=K(\mathsf v^Q)$, $\mathsf v=K^*(\mathsf u^P)$, and
\begin{equation}\label{eq:it-auxiliary-mass}
S:=\int\mathsf u^{P+1}=\int\mathsf v^{Q+1}=C_\sharp^{-n/\alpha}.
\end{equation}
Multiplying the two components by $\omega^{-1/(P+1)}$ and $\omega^{-1/(Q+1)}$, respectively, gives the coefficient-one solution for the unnormalized kernel $G$. This transformation preserves the conditional averaging operators and their probability measures, so the preceding radial spectral theorem also applies to the normalized kernel. For example, the integral kernels on Lebesgue $L^2$ agree, since
\[
\omega^{-\{a(P-1)+b(Q-1)\}/(2n)}
=\omega^{-\alpha/n}.
\]
Throughout, $S_ch(t)=h(t-c)$ denotes translation on the cylinder. If $D_\tau^{(q)}h(x)=\tau^{n/q}h(\tau x)$, then
\begin{equation}\label{eq:it-shift-radius}
\mathcal U_qD_\tau^{(q)}=S_{-\log\tau}\mathcal U_q.
\end{equation}
In particular, $S_cF$ corresponds to $D_{e^{-c}}^{(p)}f_*$: moving the cylindrical center to the right by $c$ multiplies the characteristic Euclidean radius by $e^c$. Translations chosen below to establish compactness are coordinate choices for comparing functions. The algorithm itself performs only the two norm normalizations and makes no additional spatial rescaling.

The updates below are the integral-operator analogue of the finite-dimensional nonlinear power method for mixed norms introduced by Boyd \cite{Boyd}. That precedent concerns matrices; convergence for the present translation-invariant operator requires the compactness and scale analysis that follows.

We analyze the iteration \eqref{eq:mainiteration} under the assumptions
of Theorem~\ref{thm:mainiteration}, after translating the common center
to the origin. The proof has three steps: the approximate
Euler--Lagrange equations yield a finite-profile decomposition and energy
quantization; a positive interaction estimate excludes more than one
profile; and the radial spectral gap makes the remaining translation
increments summable. The last step gives convergence to a fixed scale.

\subsection{Cylindrical formulation and monotonicity}
Use the notation of \eqref{eq:it-normalized-kernel}--\eqref{eq:it-optimal-pair}. Set $F_k=\mathcal U_p f_k$, $H_k=\mathcal U_r g_k$, and choose $(f_*,g_*)$ so that $F=\mathcal U_p f_*$ and $H=\mathcal U_r g_*$. The cylindrical functional is $\int F_kKH_k$, and hence
\begin{equation}\label{eq:iteration-cyl}
 H_k=\frac{(K^*F_k)^Q}{\|(K^*F_k)^Q\|_r},\qquad
 F_{k+1}=\frac{(KH_k)^P}{\|(KH_k)^P\|_p}
 =:\mathcal T(F_k).
\end{equation}
The extremizers satisfy \eqref{eq:it-optimal-pair}. Define the reduced energy by
\begin{equation}\label{eq:it-reduced-energy}
 E(h)=\frac{\|K^*h\|_{r'}}{\|h\|_p}\quad(h\ne0),\qquad
 E(F_k)=\mathcal J(f_k,g_k).
\end{equation}
Common dilation corresponds exactly to common translation on the cylinder:
$\mathcal U_sD_\tau^{(s)}=S_{-\log\tau}\mathcal U_s$.

\begin{lemma}[Strict increase and vanishing successive differences]\label{lem:iteration-monotone}
The iterates are nonnegative and have unit norm, and their Euclidean representatives remain radially nonincreasing about the origin. With $a_k=\|K^*F_k\|_{r'}$ and $b_k=\|KH_k\|_{p'}$, one has
\begin{equation}\label{eq:it-energy-chain}
 0<a_0\le a_k\le b_k\le a_{k+1}\le C_\sharp.
\end{equation}
Unless $(F_k,H_k)$ is an extremizing pair, $a_{k+1}>a_k$. There exists $L\in(0,C_\sharp]$ such that $a_k,b_k\to L$, and
\begin{equation}\label{eq:it-asymptotic-regularity}
 \|F_{k+1}-F_k\|_p+\|H_{k+1}-H_k\|_r\longrightarrow0.
\end{equation}
\end{lemma}
\begin{proof}
The HLS bound and strict positivity of the kernel ensure that each denominator is finite and nonzero. Convolution preserves nonnegative radially nonincreasing functions by Lieb \cite[Lemma~2.2(i)]{Lieb}; monotone convergence removes any truncations. Positive powers and division by positive constants preserve this property.

Normalization and equality in H\"older's inequality give
\begin{equation}\label{eq:it-exact-response}
 K^*F_k=a_kH_k^{r-1},\qquad KH_k=b_kF_{k+1}^{p-1}.
\end{equation}
Thus $a_k=\int F_kKH_k\le b_k$ and $b_k=\int H_kK^*F_{k+1}\le a_{k+1}$; the HLS inequality supplies the upper bound. This proves \eqref{eq:it-energy-chain} and the existence of a common positive limit.

If $a_{k+1}=a_k$, equality holds throughout these inequalities. In particular, $\int F_kF_{k+1}^{p-1}=a_k/b_k=1$, so equality in H\"older's inequality and the unit norms imply $F_k=F_{k+1}$. Equation \eqref{eq:it-exact-response} is then the normalized Euler--Lagrange system with common multiplier $a_k$. Corollary \ref{cor:stationarysharp} identifies the pair as an extremizer, proving strict increase.

In general, $\int F_kF_{k+1}^{p-1}=a_k/b_k\to1$, whence
\[
 1\ge\left\|\frac{F_k+F_{k+1}}2\right\|_p
 \ge\frac{1+a_k/b_k}{2}\longrightarrow1.
\]
Uniform convexity of $L^p$ \cite{Clarkson} gives the first vanishing difference. Similarly, $\int H_kH_{k+1}^{r-1}=b_k/a_{k+1}\to1$ gives the second. This uses only $L>0$, without presupposing convergence of the functions.
\end{proof}

\begin{lemma}[Euler--Lagrange residuals]\label{lem:it-residual}
Under the preceding assumptions,
\begin{equation}\label{eq:it-profiledefect}
 \|KH_k-LF_k^{p-1}\|_{p'}
 +\|K^*F_k-LH_k^{r-1}\|_{r'}\longrightarrow0.
\end{equation}
\end{lemma}
\begin{proof}
By \eqref{eq:it-exact-response}, the second term is $|a_k-L|$, and the first is at most
\[
 |b_k-L|+L\|F_{k+1}^{p-1}-F_k^{p-1}\|_{p'}.
\]
This tends to zero by Lemma \ref{lem:iteration-monotone} and continuity of the power map $L^p_+\to L^{p'}_+$.
\end{proof}

The critical identity \eqref{eq:it-basicparameters} gives $PQ>1$, so at least one of $P,Q$ exceeds one. When needed, interchange the components under analysis so that $P>1$; if both exceed one, also arrange $P\le Q$. This amounts to analyzing complete two-step updates of the other component of the same orbit and does not change the iteration \eqref{eq:mainiteration}. Both half-step response maps are continuous on their respective nonnegative $L^p,L^r$ spaces, so convergence transfers back to the original components. More explicitly, set
\[
\widetilde F_k=H_k,\qquad \widetilde H_k=F_{k+1},\qquad
(\widetilde p,\widetilde r,\widetilde K)=(r,p,K^*).
\]
Equation \eqref{eq:iteration-cyl} becomes a two-step iteration in precisely the same order, with initial datum $H_0$. This datum is nonnegative, has unit norm, and has a radially nonincreasing Euclidean representative. If $H_k\to S_cH$, continuity and translation equivariance of the half-step response yield $F_{k+1}\to S_cF$; shifting the index gives the same limit for $F_k$. Interchanging the components therefore neither restricts the initial datum nor changes their eventual common scale.

\section{Decomposition into separating profiles}\label{sec:iteration-compactness}

Convergence of the objective supplies only a real number $L$. To compare it with the sharp constant, one must construct limiting functions that satisfy the equations. We prove the following compactness theorem, including preservation of both unit norms.

\begin{theorem}[Strong compactness at a common scale]\label{thm:it-orbitcompact}
Under the parameter and initial-data assumptions of Theorem~\ref{thm:mainiteration}, every sequence $k_j\to\infty$ has a subsequence, still denoted $k_j$, common translation parameters $c_j\in\R$, and nonnegative functions $F_\infty\in L^p(\R)$, $H_\infty\in L^r(\R)$ such that
\begin{equation}\label{eq:it-orbitcompact}
 \|S_{-c_j}F_{k_j}-F_\infty\|_p
 +\|S_{-c_j}H_{k_j}-H_\infty\|_r\longrightarrow0,
 \qquad \|F_\infty\|_p=\|H_\infty\|_r=1.
\end{equation}
In Euclidean variables, $S_{-c_j}$ corresponds to the common norm-preserving dilation $D_{e^{c_j}}$ of the two components.
\end{theorem}

We first extract finitely many profiles from the approximate equations.
After estimating their interactions, we use square-root variables and
the energy gain to exclude more than one separating profile.

\subsection{Finite-profile decomposition of approximate stationary points}\label{sec:iteration-profiles}

The extraction of separating profiles follows the concentration--compactness strategy of Lions \cite[Lemma~I.1]{Lions}. Here the vanishing equation residuals identify every nonzero local limit as an actual Euler--Lagrange solution. Corollary \ref{cor:stationarysharp} then identifies these limits as translates and multiples of the same extremizer. Each carries a fixed positive mass, so there can be only finitely many. The residual equations and the quantization formula below require the additional verifications given in this subsection.

Retain the parameters, initial datum, and iteration of Section \ref{sec:iteration}, and write $F_k=\mathcal U_p f_k$, $H_k=\mathcal U_r g_k$, $a_k=\|K^*F_k\|_{r'}$, and $b_k=\|KH_k\|_{p'}$. Lemma \ref{lem:iteration-monotone} gives $a_k,b_k\to L\in(0,C_{\sharp}]$ and $\|F_{k+1}-F_k\|_p+\|H_{k+1}-H_k\|_r\to0$. Lemma \ref{lem:it-residual} derives \eqref{eq:it-profiledefect} from the actual iteration. The decomposition below uses these approximate equations and does not require individual iterates to solve the stationary system.

\begin{theorem}[Finite separating profiles and energy quantization]\label{thm:it-profiles}
For every initial datum satisfying Theorem~\ref{thm:mainiteration}, there is an integer $m\ge1$, uniquely determined by $L$, such that
\begin{equation}\label{eq:it-quantization}
 L=C_{\sharp}m^{-\alpha/n}.
\end{equation}
Fix a unit-norm extremizing pair $(F,H)$. Every $k_j\to\infty$ has a subsequence and real numbers $t_{1,j},\ldots,t_{m,j}$ such that $|t_{i,j}-t_{l,j}|\to\infty$ for $i\ne l$ and
\begin{equation}\label{eq:it-profilestrong}
 \left\|F_{k_j}-m^{-1/p}\sum_{i=1}^m S_{t_{i,j}}F\right\|_p
 +\left\|H_{k_j}-m^{-1/r}\sum_{i=1}^m S_{t_{i,j}}H\right\|_r
 \longrightarrow0.
\end{equation}
\end{theorem}

This decomposition has a direct radial interpretation. By \eqref{eq:it-shift-radius},
\[
S_{t_{i,j}}F
=\mathcal U_p\!\left(D_{e^{-t_{i,j}}}^{(p)}f_*\right).
\]
All profiles retain their common Euclidean center at the origin, while their characteristic radii are multiplied by $e^{t_{i,j}}$. The condition $|t_{i,j}-t_{l,j}|\to\infty$ means that the corresponding radius ratios tend to zero or infinity. More precisely, for a sufficiently large fixed $R$, each annulus
\[
e^{t_{i,j}-R}<|x|<e^{t_{i,j}+R}
\]
carries approximately $1/m$ of the $L^p$ mass; taking first $j\to\infty$ and then $R\to\infty$ makes the remaining mass arbitrarily small. We must exclude the persistent coexistence of two or more such dilation profiles, all with the same spatial center.

\begin{proof}
\textit{Step 1: uniform local regularity after updating arbitrary data.}
No derivative assumption on $F_0$ is made. Set
\[
W_k=K^*F_k,\qquad Z_k=KH_k.
\]
Lemma \ref{lem:direct-kernel-factor} gives $|\mathcal K'|\le C\mathcal K$ and $\mathcal K,\mathcal K'\in L^{p'}\cap L^{r'}$. Thus the convolutions have first weak derivatives, H\"older's inequality ensures finiteness of the relevant integrals, and
\[
|W_k'|\le CW_k,\qquad |Z_k'|\le CZ_k.
\]
Strict positivity of the kernel and nonnegative unit-norm inputs imply $W_k,Z_k>0$ everywhere. They have positive lower bounds on compact intervals, so the usual chain rule applies to the positive powers in the update. By \eqref{eq:it-exact-response},
\[
H_k=a_k^{-Q}W_k^Q,\qquad F_{k+1}=b_k^{-P}Z_k^P,
\]
and hence
\begin{equation}\label{eq:it-profiles-expanded-regularity}
|H_k'|\le CQH_k,\qquad |F_{k+1}'|\le CPF_{k+1}.
\end{equation}
The normalization factors cancel in these relative derivative estimates. The constants are therefore independent of $k$ and require no proximity of the initial energy to the sharp value.

If $h>0$, $|h'|\le C_0h$, and $\|h\|_s=1$, then, for every $t,t_0$,
\[
h(t)\ge h(t_0)e^{-C_0|t-t_0|}.
\]
Integrating its $s$th power gives $1\ge 2h(t_0)^s/(sC_0)$, so $\|h\|_\infty$ is uniformly bounded. The relative derivative bound then also bounds the first derivative. Consequently, for $k\ge1$, the functions $F_k,H_k$ and all their common translates are uniformly bounded and share a Lipschitz constant on every fixed compact interval. Arzel\`a--Ascoli and diagonal extraction give locally uniformly convergent subsequences. The unit $L^p,L^r$ bounds and reflexivity also yield weakly convergent subsequences, whose weak limits coincide with their locally uniform limits. Every local limit $(X,Y)$ after common translation therefore satisfies
\[
X,Y\ge0,\qquad \|X\|_p\le1,\qquad \|Y\|_r\le1.
\]

\textit{Step 2: every nonzero local limit carries the same positive mass.}
Suppose that, along a subsequence and after common translations, $F_{k_j},H_{k_j}$ converge locally uniformly to $X,Y$. Suppress the translations in this step. For $|t|\le A$, split the integral defining $KH_{k_j}(t)$ into $|t-z|\le R$ and $|t-z|>R$. The error in the first part is at most
\[
\|\mathcal K\|_1
\sup_{|z|\le A+R}|H_{k_j}(z)-Y(z)|,
\]
which tends to zero for fixed $R$. Applying H\"older's inequality separately to $H_{k_j}$ and $Y$ bounds the total error in the second part by
\[
2\|\mathcal K\mathbf1_{\{|t|>R\}}\|_{r'}.
\]
This is independent of $j$ and the output point and tends to zero as $R\to\infty$. Thus $KH_{k_j}\to KY$ locally uniformly. The same argument for $K^*F_{k_j}$ uses the unit $L^p$ bound and the kernel's $L^{p'}$ tail norm.

Since $p-1,r-1>0$, the power functions are uniformly continuous on bounded nonnegative intervals. Hence $F_{k_j}^{p-1}\to X^{p-1}$ and $H_{k_j}^{r-1}\to Y^{r-1}$ locally uniformly. Meanwhile, Lemma \ref{lem:it-residual} gives convergence of the residuals to zero in $L^{p'},L^{r'}$. Convolution commutes with common translations, and translations preserve these norms, so the same residual conclusion holds for the translated components. The locally uniform limit of each residual has zero norm on every compact interval; being continuous, it vanishes everywhere. Therefore
\[
KY=LX^{p-1},\qquad K^*X=LY^{r-1}.
\]
This argument uses existing locally uniform limits of the individual equation terms; it does not infer uniform convergence of residuals from their norm convergence.

If $X$ is nonzero, the second equation and positivity of the kernel give $Y>0$, and the first then gives $X>0$. The same conclusion follows if $Y$ is initially known to be nonzero. Pair the equations with $X,Y$, respectively. The HLS bound ensures finiteness, and Fubini's theorem gives
\[
L\int X^p=\int XKY=\int YK^*X=L\int Y^r.
\]
As $L>0$, denote the common positive mass by
\[
\theta:=\int X^p=\int Y^r>0.
\]
Dividing $X,Y$ by $\theta^{1/p},\theta^{1/r}$ produces unit-norm equations with common multiplier $L\theta^{1-1/p-1/r}=L\theta^{-\alpha/n}$. Corollary \ref{cor:stationarysharp} and its stationary classification give
\begin{equation}\label{eq:it-profilemass}
\theta=(L/C_\sharp)^{n/\alpha},\qquad
(X,Y)=(\theta^{1/p}S_cF,\theta^{1/r}S_cH)
\end{equation}
for some fixed $c\in\R$. In particular, $\theta>0$ depends only on the energy limit $L$ of the entire orbit, and is independent of the subsequence and local limit. We use this value of $\theta$ throughout.

\textit{Step 3: successive extraction terminates after finitely many profiles.}
We use the Brezis--Lieb lemma \cite[Theorem~1]{BrezisLieb}: if $s>1$, $h_j\to h$ almost everywhere, and $\sup_j\|h_j\|_s<\infty$, then
\begin{equation}\label{eq:it-normsplitting}
\|h_j\|_s^s-\|h_j-h\|_s^s\longrightarrow\|h\|_s^s.
\end{equation}
Fatou's lemma gives $h\in L^s$. Thus the almost-everywhere convergence and uniform $L^s$ bound are exactly the hypotheses of the cited theorem, yielding \eqref{eq:it-normsplitting}.

Along the prescribed time subsequence, suppose that $l$ common centers $t_{1,j},\ldots,t_{l,j}$ have been extracted, with pairwise distances tending to infinity, and that near each center the original components converge locally uniformly to $\theta^{1/p}F,\theta^{1/r}H$. Define
\[
R_j^{(l)}
=F_{k_j}-\theta^{1/p}\sum_{i=1}^lS_{t_{i,j}}F.
\]
For $l=0$, the sum is empty and $R_j^{(0)}=F_{k_j}$. We maintain inductively
\begin{equation}\label{eq:it-profiles-expanded-mass}
\|R_j^{(l)}\|_p^p\longrightarrow1-l\theta.
\end{equation}
For $l=0$, this follows from the unit norm.

The remainder can change sign; we first show that its negative part cannot retain a positive height. For every fixed $R>0$, on $|t-t_{i,j}|\le R$ the original function approaches the $i$th reference profile, while all other reference profiles converge uniformly to zero by separation. Thus $R_j^{(l)}$ converges uniformly to zero on the union of these intervals. Outside that union, $F_{k_j}\ge0$ gives
\[
\bigl(R_j^{(l)}(t)\bigr)_-
\le\theta^{1/p}\sum_{i=1}^lF(t-t_{i,j})
\le\theta^{1/p}l\sup_{|x|>R}F(x).
\]
The fixed extremizer $F$ tends to zero at both ends of the line. Taking first $j\to\infty$ and then $R\to\infty$ yields
\[
\|(R_j^{(l)})_-\|_\infty\longrightarrow0.
\]
Only nonnegativity of the original functions is used here; the remainder is not assumed nonnegative.

If $\|R_j^{(l)}\|_\infty\to0$, extraction stops. Otherwise, pass to a subsequence realizing a positive upper limit. Since the negative part vanishes, there are $\varepsilon_0>0$ and points $z_j$ with $R_j^{(l)}(z_j)\ge\varepsilon_0$. Their distance from every old center must tend to infinity: otherwise a subsequence would lie in a fixed neighborhood of an old center, contradicting uniform vanishing of the remainder there. Since $F_{k_j}(z_j)\ge R_j^{(l)}(z_j)$, common translation of the original components to $z_j$ and Step 1 yield a new local limit $(X,Y)$ with $X(0)\ge\varepsilon_0$. Step 2 gives the form \eqref{eq:it-profilemass}. Absorb its fixed shift $c$ into the new center by setting $t_{l+1,j}=z_j+c$. The limit near this center is exactly $(\theta^{1/p}F,\theta^{1/r}H)$; the fixed added constant preserves separation from all old centers.

At the new center, the finitely many old reference profiles tend locally uniformly to zero, so
\[
S_{-t_{l+1,j}}R_j^{(l)}
\longrightarrow\theta^{1/p}F
\quad\text{locally uniformly}.
\]
Applying \eqref{eq:it-normsplitting} gives
\[
\|R_j^{(l)}\|_p^p-\|R_j^{(l+1)}\|_p^p\longrightarrow\theta.
\]
The induction hypothesis therefore yields $\|R_j^{(l+1)}\|_p^p\to1-(l+1)\theta$, preserving \eqref{eq:it-profiles-expanded-mass}.

The limiting remainder norm cannot be negative, so at most $\lfloor1/\theta\rfloor$ profiles can be extracted. If the $L^\infty$ remainder still failed to vanish after that many extractions, the same procedure would produce an additional profile and a negative right-hand side in \eqref{eq:it-profiles-expanded-mass}, a contradiction. After finitely many subsequence extractions, there are therefore $m\ge0$ profiles such that
\begin{equation}\label{eq:it-diffuseremainder}
R_j=F_{k_j}-\theta^{1/p}\sum_{i=1}^mS_{t_{i,j}}F,\qquad
\sup_j\|R_j\|_p<\infty,\qquad \|R_j\|_\infty\to0.
\end{equation}
Neither $L^p$ vanishing nor $m\theta=1$ has yet been established. Both require the actual update equations.

\textit{Step 4: strong approximation of the second component from the first half-step.}
Since $\mathcal K\in L^1$, Young's inequality and \eqref{eq:it-diffuseremainder} give
\[
\|K^*R_j\|_p\le\|\mathcal K\|_1\|R_j\|_p\le C,\qquad
\|K^*R_j\|_\infty
\le\|\mathcal K\|_1\|R_j\|_\infty\longrightarrow0.
\]
The critical relation gives $1/r'=1/p-\alpha/n$, so $r'>p$. Thus
\[
\|K^*R_j\|_{r'}
\le\|K^*R_j\|_p^{p/r'}
   \|K^*R_j\|_\infty^{1-p/r'}
\longrightarrow0.
\]
It is the convolution remainder that vanishes in $L^{r'}$; this step does not infer $L^p$ vanishing of the remainder from its $L^\infty$ vanishing.

The exact update $K^*F_{k_j}=a_{k_j}H_{k_j}^{r-1}$ and the extremizer equation $K^*F=C_\sharp H^{r-1}$ imply
\[
\left\|H_{k_j}^{r-1}
-\frac{C_\sharp\theta^{1/p}}{a_{k_j}}
  \sum_{i=1}^mS_{t_{i,j}}(H^{r-1})\right\|_{r'}
\longrightarrow0.
\]
Since $a_{k_j}\to L$, \eqref{eq:it-profilemass} gives
\[
\frac{C_\sharp\theta^{1/p}}{a_{k_j}}
\longrightarrow
\theta^{\,1/p-\alpha/n}=\theta^{1/r'}.
\]
The $L^{r'}$ norms of the finite sums are uniformly bounded, so the coefficient may be replaced by its limit.

The positive-power continuity used here and below has the following quantitative form. For nonnegative $A,B\in L^t$, $\beta>0$, and $t/\beta>1$,
\[
\|A^\beta-B^\beta\|_{t/\beta}
\le
\begin{cases}
\|A-B\|_t^\beta,&0<\beta\le1,\\
C_\beta(\|A\|_t+\|B\|_t)^{\beta-1}\|A-B\|_t,&\beta>1.
\end{cases}
\]
The first inequality follows from the scalar H\"older estimate for powers; the second follows from the mean value theorem and H\"older's inequality. Taking $t=r'=Q+1$ and $\beta=Q$ yields
\[
H_{k_j}
-\theta^{1/r}
 \left(\sum_{i=1}^mS_{t_{i,j}}(H^{r-1})\right)^Q
\longrightarrow0
\quad\text{in }L^r.
\]

The power of the sum must still be replaced by the sum of the powers. Put $h=H^{r-1}\in L^{r'}$, and approximate $h$ in $L^{r'}$ by nonnegative compactly supported functions $h_R$. For fixed $R$, separation makes the supports of $S_{t_{i,j}}h_R$ eventually disjoint, so pointwise
\[
\left(\sum_iS_{t_{i,j}}h_R\right)^Q
=\sum_iS_{t_{i,j}}(h_R^Q).
\]
The preceding power estimate controls the replacement error uniformly in $L^r$. Taking first $j\to\infty$ and then $R\to\infty$ gives
\[
\left(\sum_iS_{t_{i,j}}(H^{r-1})\right)^Q
-\sum_iS_{t_{i,j}}H\longrightarrow0
\quad\text{in }L^r.
\]
Consequently,
\begin{equation}\label{eq:it-secondprofile}
H_{k_j}-\theta^{1/r}\sum_{i=1}^mS_{t_{i,j}}H
\longrightarrow0\quad\text{in }L^r.
\end{equation}
The argument includes $Q<1$, since all powers act on nonnegative functions. If $m=0$, \eqref{eq:it-secondprofile} would imply $H_{k_j}\to0$ in $L^r$, contradicting $\|H_{k_j}\|_r=1$. Thus $m\ge1$.

\textit{Step 5: strong approximation of the first component and determination of the number of profiles.}
Apply $K$ to \eqref{eq:it-secondprofile}. The HLS bound gives
\[
KH_{k_j}
-C_\sharp\theta^{1/r}
 \sum_{i=1}^mS_{t_{i,j}}(F^{p-1})
\longrightarrow0\quad\text{in }L^{p'}.
\]
Using the exact update $KH_{k_j}=b_{k_j}F_{k_j+1}^{p-1}$, the limit $b_{k_j}\to L$, and
\[
\frac{C_\sharp\theta^{1/r}}{L}
=\theta^{\,1/r-\alpha/n}=\theta^{1/p'},
\]
we obtain
\[
F_{k_j+1}^{p-1}
-\theta^{1/p'}\sum_{i=1}^mS_{t_{i,j}}(F^{p-1})
\longrightarrow0\quad\text{in }L^{p'}.
\]
Apply the preceding power estimate with $t=p'=P+1$, $\beta=P$, and the same compact-support approximation for the finite separating translates. This gives
\[
F_{k_j+1}
-\theta^{1/p}\sum_{i=1}^mS_{t_{i,j}}F
\longrightarrow0\quad\text{in }L^p.
\]
Lemma \ref{lem:iteration-monotone} independently gives $\|F_{k_j+1}-F_{k_j}\|_p\to0$, so this also holds for $F_{k_j}$, that is,
\[
\|R_j\|_p\longrightarrow0.
\]
Strong approximation of both components has now been obtained from their actual half-step equations.

The $p$th powers of the norms of finitely many separating translates add asymptotically. This follows again by nonnegative compact-support approximation or repeated application of \eqref{eq:it-normsplitting}. Hence
\[
\left\|\theta^{1/p}\sum_{i=1}^mS_{t_{i,j}}F\right\|_p^p
\longrightarrow m\theta.
\]
Strong $L^p$ approximation and $\|F_{k_j}\|_p=1$ force the left-hand side to tend to one. Thus $m\theta=1$. Substitution into \eqref{eq:it-profilemass} gives $L=C_\sharp m^{-\alpha/n}$, $\theta^{1/p}=m^{-1/p}$, and $\theta^{1/r}=m^{-1/r}$, proving \eqref{eq:it-quantization} and \eqref{eq:it-profilestrong}. Since $L$ is the fixed energy limit of the entire orbit, it determines $m$ uniquely, independently of the time subsequence.
\end{proof}

\begin{corollary}[Reference parameters for the entire sequence]\label{cor:it-centers}
Let $m$ be the integer in Theorem \ref{thm:it-profiles}. If $m=1$, then $\inf_{c\in\R}(\|F_k-S_cF\|_p+\|H_k-S_cH\|_r)\to0$. If $m\ge2$, then, for every sufficiently large $k$, one may choose $t_{1,k}<\cdots<t_{m,k}$ such that
\begin{equation}\label{eq:it-orderedcenters}
 D(\mathbf t_k):=\min_{i\ne l}|t_{i,k}-t_{l,k}|\to\infty,
 \qquad \|F_k-F_{\mathbf t_k}\|_p\to0,
 \qquad \max_i|t_{i,k+1}-t_{i,k}|\to0,
\end{equation}
where $F_{\mathbf t}=\sum_iS_{t_i}F/\|\sum_iS_{t_i}F\|_p$.
\end{corollary}
\begin{proof}
The energy limit $L$ is the same for the entire orbit, so every subsequence has the same number of reference profiles. For each fixed $R>0$, the $L^p\times L^r$ distance to the family of equal-weight $m$-profile sums with separation at least $R$ tends to zero. Otherwise, a subsequence with distance bounded below would contradict Theorem \ref{thm:it-profiles}. Choose a stepwise increasing $R_k\to\infty$ and ordered centers in the corresponding family whose distance is at most the infimum plus $1/k$. Norm splitting for finite separating translates makes the normalization coefficient tend to $m^{-1/p}$, giving the first two assertions in \eqref{eq:it-orderedcenters}. For $m=1$, the distance conclusion directly gives the asserted vanishing infimum; convergence of the chosen translation parameters is not yet required.

Suppose $m\ge2$. Asymptotic regularity and the first two assertions imply that consecutive reference functions approach one another in $L^p$. Translate around any center of the first group. If the second group had no center within bounded distance, it would tend to zero on every fixed compact interval, whereas the first would tend to $m^{-1/p}F$, a contradiction. There is therefore a matching center at bounded distance, and separation makes the match unique. If the distance tends along a subsequence to $\delta$, the local limits give $F=S_\delta F$. A nonzero $L^p(\R)$ function cannot have a nonzero period, so $\delta=0$. Match the finitely many centers in this way. Ordering from left to right preserves their indices, so the maximum center displacement tends to zero.
\end{proof}

\subsection{Energy and update errors of separating profile sums}\label{sec:iteration-interaction}

The finite-profile decomposition alone does not exclude $m\ge2$, since sufficiently separated profiles can still approximately solve the equations. We use the direction of the actual iteration: its nondecreasing energy approaches its limit from below. We first prove that the normalized sum of $m$ separating profiles has energy $L_m+c_m\varepsilon+o(\varepsilon)$ with $c_m,\varepsilon>0$. Its update error is smaller than this positive correction in the sense needed for the energy estimates below.

The reference family consists of finite sums of translates of a fixed extremizer and is not itself an iteration orbit. Section~\ref{sec:root-method} uses square-root variables to bound the difference of the two energies by the actual one-step energy gain and the squared reference update error. Combined with the estimates here, this forces the nonnegative energy deficit to decrease faster than the slowly varying interaction permits. A scalar recurrence then excludes $m\ge2$ and proves Theorem~\ref{thm:it-orbitcompact}.

Fix $n\ge3$, $0<\lambda\le n-2$, and a finite number $m\ge2$ of reference profiles. To estimate the full two-step map, choose the component with $P>1$; if both powers exceed one, also arrange $P\le Q$. As explained above, this merely analyzes complete two-step updates of the other component, interchanging $K,K^*$ if necessary, and leaves the algorithm unchanged. The conclusion for the other component follows by half-step continuity. Thus $1<p<2$, and the critical relations give
\begin{equation}\label{eq:it-parameteridentities}
 a+b=\lambda,\qquad aP=b+\alpha,\qquad bQ=a+\alpha,
 \qquad \alpha=n-\lambda\ge2.
\end{equation}
Use the auxiliary profiles $(\mathsf u,\mathsf v)$ from \eqref{eq:it-auxiliary-profiles}. They satisfy the coefficient-one equations, and multiplication by fixed positive amplitude constants preserves tail rates and positive-coefficient asymptotics. Retain the notation
\[
 S=\int_\R \mathsf u^{P+1}\dd t=\int_\R \mathsf v^{Q+1}\dd t,
 \qquad F=S^{-1/p}\mathsf u^P,\qquad H=S^{-1/r}\mathsf v^Q.
\]
The amplitude conversion gives $C_{\sharp}=S^{-\alpha/n}$. Write $\mathsf u_i=S_{t_i}\mathsf u$, $\mathsf v_i=S_{t_i}\mathsf v$, $F_i=S_{t_i}F$, $H_i=S_{t_i}H$, and set
\begin{align}
 D(\mathbf t)&=\min_{i\ne j}|t_i-t_j|,\qquad
 F_{\mathbf t}=\frac{\sum_iF_i}{\|\sum_iF_i\|_p},\label{eq:it-reference}\\
 \varepsilon(\mathbf t)&=\int_\R\left[(\sum_i \mathsf v_i)^{Q+1}-\sum_i \mathsf v_i^{Q+1}\right]\dd t,\label{eq:it-interactionquantities}\\
 \zeta(\mathbf t)&=\|\mathcal T(F_{\mathbf t})-F_{\mathbf t}\|_p.\notag
\end{align}
All uniformity statements keep $m$ and the parameters fixed while allowing arbitrary center configurations.

\subsubsection{Positive energy correction from pairwise overlap}

\begin{lemma}[Overlap integrals for two-sided tails]\label{lem:it-overlap}
Let $w>0$ be continuous, with left and right logarithmic tail rates $\ell,\rho>0$. For $q>1$, define
\[
 R_q(x,y)=(x+y)^q-x^q-y^q,\qquad
 J_q^w(D)=\int_\R R_q(w(t),w(t-D))\dd t.
\]
Then $J_q^w(D)>0$ and
\begin{equation}\label{eq:it-overlaprate}
 \lim_{D\to\infty}-D^{-1}\log J_q^w(D)
 =\sigma_q(\ell,\rho)
 :=\min\left\{\ell,\rho,\frac{q\ell\rho}{\ell+\rho}\right\}.
\end{equation}
If both tails admit two-sided exponential comparisons, with fixed-degree polynomial factors allowed in the upper bounds, then, for $c>0$,
\begin{equation}\label{eq:it-overlappower}
 \int_\R R_q(w,w_D)^c\dd t
 \le C(1+D)^M\exp[-c\sigma_q(\ell,\rho)D].
\end{equation}
For a fixed number of terms and $q>1$,
\begin{equation}\label{eq:it-finitepaircomparison}
 (\sum_i x_i)^q-\sum_i x_i^q
 \asymp_{m,q}\sum_{i<j}R_q(x_i,x_j),\qquad x_i\ge0.
\end{equation}
\end{lemma}
\begin{proof}
For $x\ge y>0$, write the quotient $R_q(x,y)/(x^{q-1}y)$ as $((1+s)^q-1-s^q)/s$, with $0<s\le1$. It follows that
\[
 R_q(x,y)\asymp_q\min(x,y)\max(x,y)^{q-1}.
\]
First consider the model $w(t)=e^{\ell t}$ for $t\le0$ and $w(t)=e^{-\rho t}$ for $t\ge0$. On $[0,D]$, the two tails cross at $t_0=\ell D/(\ell+\rho)$. On the two sides of the crossing, the negative logarithms of the overlap are respectively
\[
 \ell D+[(q-1)\rho-\ell]t,\qquad
 (q-1)\ell D+[\rho-(q-1)\ell]t.
\]
Their endpoint values are $\ell D,q\ell\rho D/(\ell+\rho),\rho D$. The exterior half-lines contribute $e^{-\ell D}$ and $e^{-\rho D}$, respectively; a zero slope adds only a factor $D$. For general logarithmic tail rates, sandwich between models with rates $\ell\pm\delta,\rho\pm\delta$ and let $\delta\downarrow0$. The same four-region integration after taking the $c$th power proves \eqref{eq:it-overlappower}. Polynomial tail factors remain of bounded degree under this decomposition and are absorbed into $(1+D)^M$.

Finally, the middle expression in \eqref{eq:it-finitepaircomparison} is at least each pairwise interaction, giving the lower bound. Choose a maximal $x_k$. The mean value theorem bounds that expression above by $q(mx_k)^{q-1}\sum_{j\ne k}x_j$; the two-variable comparison then proves the upper bound.
\end{proof}

\begin{proposition}[Positive interaction for finite reference sums]\label{prop:it-interaction}
Let $P>1$ and fix $m\ge2$. With $L_m=C_{\sharp}m^{-\alpha/n}$, as $D(\mathbf t)\to\infty$,
\begin{equation}\label{eq:it-positiveinteraction}
 E(F_{\mathbf t})
 =L_m\left[1+\frac{\varepsilon(\mathbf t)}{mS(Q+1)}
       +o(\varepsilon(\mathbf t))\right]>L_m,
 \qquad E(h)=\frac{\|K^*h\|_{r'}}{\|h\|_p}.
\end{equation}
The remainder is uniform over the entire separated reference family.
\end{proposition}
\begin{proof}
By the solution-tail estimates in Proposition \ref{prop:direct-actual-tails}, the left tail rates of $\mathsf u,\mathsf v$ are $a,b$, and the right tail rates are
\[
 \rho_u=\min\{b,aQ\},\qquad \rho_v=\min\{a,bP\}.
\]
Resonant polynomial factors do not alter logarithmic rates. The overlap rates for the $p$th power of the denominator profile $\mathsf u^P$ and the $(Q+1)$st power of the numerator profile $\mathsf v$ are, respectively,
\[
 s_f=\min\left\{Pa,P\rho_u,
          \frac{(P+1)a\rho_u}{a+\rho_u}\right\},
 \qquad s_v=\sigma_{Q+1}(b,\rho_v)\le\min(a,b).
\]
Each term defining $s_f$ exceeds $s_v$. Indeed, $Pa=b+\alpha>b$ and $P\rho_u=\min(Pb,PQa)>\min(a,b)$. The third term exceeds $b$ when $\rho_u=b$, since $Pa>b$, and exceeds $a$ when $\rho_u=aQ$, since $PQ>1$. Thus $J_p^{\mathsf u^P}(D)=o(J_{Q+1}^{\mathsf v}(D))$, uniformly for all sufficiently large $D$ because of the strict exponential gap. Set
\[
 \Delta_f=\int_\R\left[(\sum_i \mathsf u_i^P)^p-\sum_i \mathsf u_i^{P+1}\right]\dd t.
\]
Equation \eqref{eq:it-finitepaircomparison} gives $\Delta_f=o(\varepsilon)$, with both quantities tending to zero and $\varepsilon>0$. The linear equation $K^*(\sum_i \mathsf u_i^P)=\sum_i \mathsf v_i$ yields the exact identity
\[
 E(F_{\mathbf t})
 =\frac{(mS+\varepsilon)^{1/(Q+1)}}{(mS+\Delta_f)^{1/p}}.
\]
Expanding the two scalar powers proves \eqref{eq:it-positiveinteraction}.
\end{proof}

\subsubsection{The update defect is smaller than the interaction energy}

\begin{proposition}[Superlinear update defect of the reference family]\label{prop:it-forcing}
Choose the component as specified at the beginning of this subsection. There exist $\eta>0$, $C<\infty$, and $D_0<\infty$ such that
\begin{equation}\label{eq:it-superlinearforcing}
 \zeta(\mathbf t)^p\le C\varepsilon(\mathbf t)^{1+\eta},
 \qquad D(\mathbf t)\ge D_0.
\end{equation}
Moreover, there is $C_{\log}<\infty$ such that any two configurations of centers satisfy
\begin{equation}\label{eq:it-loglipschitzinteraction}
 |\log\varepsilon(\mathbf t')-\log\varepsilon(\mathbf t)|
 \le C_{\log}\max_i|t_i'-t_i|.
\end{equation}
\end{proposition}
\begin{proof}
Write $Y=\sum_i \mathsf u_i^P$, $V=\sum_i \mathsf v_i$, and $W=K(V^Q)$.
All positive normalization factors at the intermediate half-step cancel exactly in the final normalization. Thus
\begin{equation}\label{eq:it-forcingnormalization}
 \mathcal T(F_{\mathbf t})=\frac{W^P}{\|W^P\|_p},\qquad
 F_{\mathbf t}=\frac{Y}{\|Y\|_p},\qquad
 \zeta\le\frac{2\|W^P-Y\|_p}{\min(\|W^P\|_p,\|Y\|_p)}.
\end{equation}
Since $W\ge \mathsf u_i$ and $Y\ge \mathsf u_i^P$, the denominator has a fixed positive lower bound. We retain estimates for the $p$th power of the unnormalized error.

We first record the precise tail bounds, which follow directly from Proposition~\ref{prop:direct-actual-tails} at an extremizer.
If $0<Q\le1<P$, then $b>a$ and
\begin{equation}\label{eq:it-mixedtails}
 \mathsf u(t)\asymp\begin{cases}e^{at},&t\le0,\\e^{-aQt},&t\ge0,\end{cases}
 \qquad
 \mathsf v(t)\asymp\begin{cases}e^{bt},&t\le0,\\e^{-at},&t\ge0.\end{cases}
\end{equation}
If $1<P\le Q$, then $a\ge b$; the two tail rates of $\mathsf u$ are $a,b$, and those of $\mathsf v$ are $b,\rho_v$, where $\rho_v=\min(a,Pb)\ge b$.
Only when $a=Pb$ does the right-tail upper bound for $\mathsf v$ contain a factor $1+t$; this factor is included in $(1+D)^M$ in all estimates below.
The piecewise integrations in Lemma~\ref{lem:it-overlap} yield
\begin{equation}\label{eq:it-energykappa}
 J_{Q+1}^{\mathsf v}(D)\asymp e^{-\kappa D},\qquad
 \kappa=\begin{cases}a,&0<Q\le1<P,\\b,&1<P\le Q.\end{cases}
\end{equation}
The smallest rate occurs at a nonresonant endpoint: in the first case $(Q+1)ab/(a+b)>a$, and in the second $(Q+1)b\rho_v/(b+\rho_v)>b$.
The resonant right-tail rate is $\rho_v=Pb>b$, so its polynomial factor does not affect this two-sided comparison.

\emph{Two-profile estimate for mixed powers $0<Q<1<P$.}
Place the centers at $0,D$, write $Y=\mathsf u^P+\mathsf u_D^P$ and $W=K((\mathsf v+\mathsf v_D)^Q)$, and set
\[
 t_u=\frac{D}{Q+1},\qquad t_v=\frac{bD}{a+b},\qquad t_u<t_v.
\]
The tail estimates give $\mathsf u_D\le C\mathsf u$ for $t\le t_u$ and $\mathsf u\le C\mathsf u_D$ for $t\ge t_u$.
To retain the cancellation between the inner concave power and the slow tail of the right-hand reference profile, write
\[
 \delta_R=W-\mathsf u_D=Kq_R,\qquad
 0\le q_R=(\mathsf v+\mathsf v_D)^Q-\mathsf v_D^Q
       \le\min(\mathsf v^Q,Q\mathsf v_D^{Q-1}\mathsf v).
\]
Substitution of \eqref{eq:it-mixedtails} gives
\begin{equation}\label{eq:it-mixedsource}
 q_R(s)\le C\begin{cases}
 e^{Qbs},&s\le0,\\
 e^{-aQs},&0\le s\le t_v,\\
 e^{-as+(1-Q)b(D-s)},&t_v\le s\le D,\\
 e^{-aD-aQ(s-D)},&s\ge D.
 \end{cases}
\end{equation}
Lemma~\ref{lem:direct-kernel-factor} gives $\mathcal K(x)\asymp e^{ax}$ for $x\le0$ and $\mathcal K(x)\asymp e^{-bx}$ for $x\ge0$.
Splitting the convolution into the four regions in \eqref{eq:it-mixedsource} and at $s=t$ yields
\begin{equation}\label{eq:it-mixedpotential}
 \delta_R(t)\le C(1+D)\begin{cases}
 e^{-aQt},&t_u\le t\le t_v,\\
 e^{-at+(1-Q)b(D-t)},&t_v\le t\le D,\\
 e^{-aD-aQ(t-D)},&t\ge D.
 \end{cases}
\end{equation}
In the first region, the leading contribution comes from $s=t$ within $[0,t_v]$; the source to the left contributes $e^{-bt}\le e^{-aQt}$.
In the second region, the preceding source decays at rate $b$ and the current source at rate $a+(1-Q)b$; their difference is $Qb-a=\alpha>0$.
In the third region, the source to the left still decays at rate $b$, while the tail source decays at rate $aQ<b$.
A zero slope in any of these integrals contributes at most the displayed polynomial factor; these computations also control all contributions across the regions.
Concavity further gives $0\le W-\mathsf u_D\le \mathsf u$ and $0\le W-\mathsf u\le \mathsf u_D$.
The mean value theorem therefore implies
\begin{equation}\label{eq:it-mixedpointwiseerror}
 |W^P-Y|\le C\begin{cases}
 \mathsf u^{P-1}\mathsf u_D+\mathsf u_D^P,&t\le t_u,\\
 \mathsf u_D^{P-1}\delta_R+\mathsf u^P,&t\ge t_u.
 \end{cases}
\end{equation}
To exhibit the strict gap in the decay rates, raise the right-hand side to the $p$th power and integrate on each region.
The left-hand cross term has rate $ap$ on $(-\infty,0]$, and endpoint rates $ap,a+\alpha$ on $[0,t_u]$.
The rates of $\mathsf u_D^{P+1}$ on $(-\infty,t_u]$ and of $\mathsf u^{P+1}$ on $[t_u,\infty)$ are both $a+\alpha$.
For the right-hand cross term, the negative logarithms on $[t_u,t_v]$ and $[t_v,D]$, after division by $p$, are respectively
\[
 a[(P-1)(D-t)+Qt],\qquad at+2\alpha(D-t),
\]
where $(P-1)a-(1-Q)b=2\alpha$.
The successive endpoint rates are $a+\alpha$, $pa(b+2\alpha)/(a+b)$, and $ap$; the interval $[D,\infty)$ again gives $ap$.
All the exponents are affine, so their endpoint values control the full integrals. Hence
\begin{equation}\label{eq:it-mixedsigma}
 \|W^P-Y\|_p^p\le C(1+D)^M e^{-\sigma D},\qquad
 \sigma=\min\left\{ap,a+\alpha,\frac{pa(b+2\alpha)}{a+b}\right\}>a.
\end{equation}
The last strict inequality follows directly from
\[
 \frac{p(b+2\alpha)}{a+b}-1
 =\frac{\alpha(a+2b+2\alpha)}{(a+b)(b+\alpha)}>0.
\]

\emph{Two-profile estimates for two superlinear powers and for a linear inner power.}
If $1<P\le Q$, write $U=\mathsf u+\mathsf u_D$ and $Z=K R_Q(\mathsf v,\mathsf v_D)$, so that $W=U+Z$.
Lemma~\ref{lem:it-overlap} gives
\begin{align}
 \|R_Q(\mathsf v,\mathsf v_D)\|_r^r&\le C(1+D)^M e^{-\sigma_gD},\qquad
 \|R_P(\mathsf u,\mathsf u_D)\|_p^p\le C(1+D)^M e^{-\sigma_fD},\notag\\
 \sigma_g&=\min\left\{rb,r\rho_v,
           \frac{(Q+1)b\rho_v}{b+\rho_v}\right\}>b,\label{eq:it-superlinearsigmag}\\
 \sigma_f&=\min\left\{pa,pb,
           \frac{(P+1)ab}{a+b}\right\}>b.\label{eq:it-superlinearsigmaf}
\end{align}
Strictness for the third terms follows respectively from $Q\rho_v>b$ and $Pa>b$.
HLS boundedness and $Pp=p'=P+1$ imply
\[
 \|Z\|_{p'}\le C\|R_Q(\mathsf v,\mathsf v_D)\|_r,\qquad
 \|(U+Z)^P-U^P\|_p
 \le C(\|U\|_{p'}^{P-1}\|Z\|_{p'}+\|Z\|_{p'}^P).
\]
Combining this with $W^P-Y=[(U+Z)^P-U^P]+R_P(\mathsf u,\mathsf u_D)$ and using $p\ge r$, we obtain
\[
 \|W^P-Y\|_p^p\le C(1+D)^M e^{-\sigma D},
 \qquad \sigma=\min(\sigma_gp/r,\sigma_f)>b.
\]
If $Q=1<P$, then $W=\mathsf u+\mathsf u_D$, and both tail rates of $\mathsf u$ equal $a$. Thus
\[
 \|W^P-Y\|_p^p\le C(1+D)^M e^{-\sigma D},
 \qquad \sigma=\min\left\{pa,\frac{(P+1)a}{2}\right\}>a.
\]

\emph{Extension to a fixed finite number of profiles.}
For mixed powers, partition the line into dominance sets $E_i=\{\mathsf u_i=\max_j \mathsf u_j\}$, resolving ties arbitrarily, and put
\[
 \delta_{ij,i}=K[(\mathsf v_i+\mathsf v_j)^Q-\mathsf v_i^Q].
\]
The increment of a concave power decreases as the base point increases. Consequently,
\[
 0\le W-\mathsf u_i\le\sum_{j\ne i}\delta_{ij,i},\qquad
 W\le\sum_j \mathsf u_j\le m \mathsf u_i\quad\hbox{on }E_i.
\]
The mean value theorem gives
\[
 |W^P-Y|\le C_m\sum_{j\ne i}(\mathsf u_i^{P-1}\delta_{ij,i}+\mathsf u_j^P)
 \quad\hbox{on }E_i.
\]
For every pair, $E_i\subset\{\mathsf u_i\ge \mathsf u_j\}$.
The two-sided tail comparisons place each actual two-profile dominance set on the corresponding side of the model interface $t_u$, up to a strip of uniformly bounded width; no monotonicity of translation ratios is required.
Integrating separately over this strip, the same endpoint rates as in \eqref{eq:it-mixedpotential}--\eqref{eq:it-mixedsigma} give
\[
 \int_{\{\mathsf u_i\ge \mathsf u_j\}}
 [(\mathsf u_i^{P-1}\delta_{ij,i})^p+\mathsf u_j^{P+1}]\dd t
 \le C(1+D_{ij})^M e^{-\sigma D_{ij}},\qquad D_{ij}=|t_i-t_j|.
\]
Summing over the finitely many ordered pairs yields the sum of the same two-profile bounds.
For two superlinear powers, apply \eqref{eq:it-finitepaircomparison} first to the inner power and then to the outer power, together with the preceding HLS and power-difference estimates. This gives
\begin{equation}\label{eq:it-finiteforcingpairs}
 \|W^P-Y\|_p^p\le C_m\sum_{i<j}(1+D_{ij})^M e^{-\sigma D_{ij}},
 \qquad \sigma>\kappa.
\end{equation}
For a linear inner power, only the outer comparison is needed.
Choose $0<\eta<\sigma/\kappa-1$, absorb the polynomial factors using \eqref{eq:it-energykappa}, and observe that
\[
 \sum_{i<j}J_{Q+1}^{\mathsf v}(D_{ij})^{1+\eta}
 \le\left(\sum_{i<j}J_{Q+1}^{\mathsf v}(D_{ij})\right)^{1+\eta}
 \asymp_m\varepsilon(\mathbf t)^{1+\eta}.
\]
Together with \eqref{eq:it-forcingnormalization}, this proves \eqref{eq:it-superlinearforcing}.

Finally, the logarithmic derivative bound for the kernel and $\mathsf v=K^*(\mathsf u^P)$ give $|\mathsf v'|\le L\mathsf v$.
Set $\delta=\max_i|t_i'-t_i|$. Each new reference term lies between $e^{-L\delta}$ and $e^{L\delta}$ times its old counterpart.
Set $C_{\log}=(Q+1)L$. The interaction function $(\sum_i x_i)^{Q+1}-\sum_i x_i^{Q+1}$ is increasing in every variable and homogeneous of degree $Q+1$.
Pointwise comparison followed by integration yields
$e^{-(Q+1)L\delta}\varepsilon(\mathbf t)\le\varepsilon(\mathbf t')\le e^{(Q+1)L\delta}\varepsilon(\mathbf t)$, proving \eqref{eq:it-loglipschitzinteraction}.
\end{proof}

\section{Square-root variables and exclusion of multiple profiles}\label{sec:root-method}

The change of variables $x=f^{p/2}$ turns the mass constraint into
$\|x\|_2=1$. The argument below has a simple structure. One update gains
energy by at least a constant times the squared change in $x$.
After adjusting the profile translations, the same change controls
the distance to the reference family. A separated reference has energy
strictly above the proposed limiting value, whereas the actual orbit
approaches that value from below. These facts force the energy deficit
to decrease faster than the slowly varying profile interaction permits.
We first prove the estimates needed for this contradiction.

\subsection{A uniform expansion in square-root variables}

The critical HLS gain permits a first-order expansion on the nonnegative
unit sphere of $L^2$, even when an individual power is below one.
The following lemma isolates the required argument.

\begin{lemma}[Expansion of a normalized positive map]\label{lem:root-halfmap}
Let $\mu,\nu$ be probability measures, $1<s<t<\infty$, and let $A$ be a
positive integral operator such that $A1=1$ and
\[
 \|A\|_{L^1(\mu)\to L^1(\nu)}
 +\|A\|_{L^s(\mu)\to L^t(\nu)}\le M.
\]
For $x\ge0$ near $1$ in $L^2(\mu)$, set
\[
 \Phi(x)=\bigl[A(x^{2/s})\bigr]^{t/2},\qquad
 \mathcal M(x)=\frac{\Phi(x)}{\|\Phi(x)\|_{L^2(\nu)}}.
\]
Then, for $1+h\ge0$ and $\delta=\|h\|_{L^2(\mu)}\to0$,
\begin{equation}\label{eq:root-halfmap}
 \mathcal M(1+h)
 =1+\frac ts\left(Ah-\int Ah\dd\nu\right)+o_{L^2(\nu)}(\delta).
\end{equation}
The remainder is uniform over all such operators and probability measures
with the same $s,t,M$. This is an expansion for nonnegative inputs; no
open neighborhood contained in the positive cone is assumed.
\end{lemma}

\begin{proof}
Positivity and $A1=1$ give $\|A\|_{\infty\to\infty}=1$.
We need a bound $A:L^q(\mu)\to L^2(\nu)$ with $1<q<2$.
These standard interpolation bounds \cite[Chapter~5, \S1]{BerghLofstrom}
follow here directly from H\"older's inequality for the nonnegative kernel.
If $t\le2$, put $k=2/t$ and $q=ks$. Since
$|Ah|^k\le A(|h|^k)$,
\[
 \|Ah\|_2\le\|A(|h|^k)\|_t^{1/k}\le C\|h\|_q.
\]
If $t>2$, put $\vartheta=t/[2(t-1)]$ and
$1/q=1-\vartheta+\vartheta/s$. For $h,v\ge0$ and $\|v\|_2=1$,
H\"older's inequality on the kernel measure gives
\[
 \langle Ah,v\rangle_\nu
 \le \langle A(h^q),1\rangle_\nu^{1-\vartheta}
     \langle A(h^{q/s}),v^{2/t'}\rangle_\nu^\vartheta
 \le C\|h\|_q.
\]
Taking the supremum over $v$ proves the bound, with
\[
 \frac2q=\kappa:=\frac{t-2+t/s}{t-1}>1.
\]
In both cases $q>1$, all constants are uniform, and the probability
normalization also gives $\|A\|_{2\to2}\le C$.

Fix $0<\epsilon<1/2$, set
$B=\{|h|>\epsilon\}$, $h_g=h\mathbf1_{B^c}$ and
$h_b=h\mathbf1_B$. Since $\mu(B)\le\delta^2/\epsilon^2$,
\begin{equation}\label{eq:root-bad-linear}
 \|Ah_b\|_2\le C\|h_b\|_q
 \le C_\epsilon\delta^{2/q}=o(\delta).
\end{equation}
The two powers can be expanded uniformly at the bounded input $1+h_g$.
Using $A1=1$ and the $L^2$ bound gives
\[
 \Phi(1+h_g)=1+\frac ts Ah_g+O_{L^2}(\epsilon\delta).
\]
It remains to control the nonlinear contribution from $B$.
Put $a=2/s$ and
$\Delta=(1+h)^a-(1+h_g)^a$.
For $z\ge-1$ and $|z|>\epsilon$, the scalar bounds
\[
 |(1+z)^a-1|\le C_\epsilon |z|^a,\qquad
 |(1+z)^a-1|\le C_\epsilon |z|^2
\]
hold because $0<a<2$. Consequently,
\[
 \|\Delta\|_s\le C_\epsilon\delta^{2/s},\qquad
 \|\Delta\|_1\le C_\epsilon\delta^2.
\]
If $t\le2$, the scalar power $t/2$ is H\"older continuous, so
\[
 \|\Phi(1+h)-\Phi(1+h_g)\|_2
 \le\|A\Delta\|_t^{t/2}
 \le C_\epsilon\delta^{t/s}=o(\delta).
\]
If $t>2$, interpolation between the $L^1$ and $L^t$ bounds gives
$\|A\Delta\|_2\le C_\epsilon\delta^\kappa$.
Since $A((1+h_g)^a)$ is uniformly bounded, the elementary power estimate
for exponent $t/2>1$ yields
\[
 \|\Phi(1+h)-\Phi(1+h_g)\|_2
 \le C\bigl(\|A\Delta\|_2+\|A\Delta\|_t^{t/2}\bigr)
 \le C_\epsilon(\delta^\kappa+\delta^{t/s})=o(\delta).
\]
Both exponents exceed one. Combining these estimates with
\eqref{eq:root-bad-linear}, then letting $\epsilon\downarrow0$, proves
$\Phi(1+h)=1+(t/s)Ah+o(\delta)$ uniformly.
The differential of $z\mapsto z/\|z\|_2$ at the unit vector $1$ is
$w\mapsto w-(\int w\dd\nu)1$, which gives \eqref{eq:root-halfmap}.
\end{proof}

\subsection{The reference operators and a static residual estimate}

Keep $1<p<2$ and put $d=PQ>1$. The square-root variables and full update are
\[
 X=F^{p/2},\qquad Y=H^{r/2},\qquad
 \mathcal N(x)=\bigl[\mathcal T(x^{2/p})\bigr]^{p/2},\qquad
 U_{\mathbf t}=F_{\mathbf t}^{p/2}.
\]
Thus $X,Y,U_{\mathbf t}$ have unit $L^2$ norm. All operator norms in this
subsection are on real $L^2$, unless other spaces are specified.
For a unit vector $Z$, write $\Pi_Z=I-Z\otimes Z$, where
$(Z\otimes Z)e=Z\langle Z,e\rangle$.

\begin{lemma}[Uniform expansion of the full root update]\label{lem:root-expansion}
Fix $m$. For every reference $f=F_{\mathbf t}$, write
\[
 g=\frac{(K^*f)^Q}{\|(K^*f)^Q\|_r},\qquad f_+=\mathcal T(f),\qquad
 c_f=\|K^*f\|_{r'},\qquad b_f=\|Kg\|_{p'},
\]
and set $X_f=f^{p/2}$, $Y_f=g^{r/2}$, and $X_+=f_+^{p/2}$.
There is a uniformly bounded operator $\mathcal L_{\mathbf t}$ such that
\begin{equation}\label{eq:root-full-expansion}
 \mathcal N(U_{\mathbf t}+e)
 =\mathcal N(U_{\mathbf t})+\mathcal L_{\mathbf t}e
       +o_{L^2}(\|e\|_2),\qquad U_{\mathbf t}+e\ge0.
\end{equation}
The remainder is uniform over all center configurations with this fixed $m$.
More precisely, define operators between the probability spaces
$\dd\mu=f^p\dd t$, $\dd\nu=g^r\dd t$, and $\dd\mu_+=f_+^p\dd t$ by
\[
 A_fh=\frac{K^*(fh)}{c_fg^{r-1}},\qquad
 D_fj=\frac{K(gj)}{b_ff_+^{p-1}}.
\]
Their conjugates to Lebesgue $L^2$ have kernels
\begin{equation}\label{eq:root-kernel-legs}
 \widehat A_f(z,t)=c_f^{-1}\mathcal K(t-z)
               f(t)^{1-p/2}g(z)^{1-r/2},\qquad
 \widehat D_f(t,z)=b_f^{-1}\mathcal K(t-z)
               f_+(t)^{1-p/2}g(z)^{1-r/2},
\end{equation}
and the principal part in \eqref{eq:root-full-expansion} is
\begin{equation}\label{eq:root-full-derivative}
 \mathcal L_{\mathbf t}=d\Pi_{X_+}\widehat D_f\widehat A_f.
\end{equation}
The same conclusions hold at the single-profile reference $f=F$.
\end{lemma}

\begin{proof}
The finite-sum comparison
\[
 m^{\min(q-1,0)}\sum_i a_i^q
 \le\Big(\sum_i a_i\Big)^q
 \le m^{\max(q-1,0)}\sum_i a_i^q\qquad(q>0)
\]
applied first with $q=Q$ and then with $q=P$ gives
\[
 R_m^{-1}f\le f_+\le R_m f,\qquad
 R_m=m^{P|Q-1|+P-1}.
\]
Also $f\ge F_i/m$, so $c_f\ge C_\sharp/m$, and the best-response
identities imply $b_f\ge c_f$. Both upper bounds are $C_\sharp$.
The HLS inequality gives uniform bounds
\[
 A_f:L^p(\mu)\longrightarrow L^{r'}(\nu),\qquad
 D_f:L^r(\nu)\longrightarrow L^{p'}(\mu_+).
\]
Both operators are positive and preserve constants. Fubini's theorem gives
\[
 \|A_f\|_{1\to1}
 \le\frac{b_f}{c_f}\sup(f_+/f)^{p-1},\qquad
 \|D_f\|_{1\to1}
 \le\frac{c_{f_+}}{b_f}\sup(g_+/g)^{r-1},
\]
where $g_+$ is the normalized response to $f_+$.
Comparability of $f_+$ and $f$ gives
$R_m^{-2Q}g\le g_+\le R_m^{2Q}g$, so these bounds are uniform.
The critical relation gives $r'>p$ and $p'>r$.
Lemma~\ref{lem:root-halfmap} therefore applies to $A_f$ with
$(s,t)=(p,r')$ and to $D_f$ with $(s,t)=(r,p')$.
After conjugation by the square roots of the measures, the two principal
parts are respectively
\[
 \frac{r'}p\Pi_{Y_f}\widehat A_f,
 \qquad \frac{p'}r\Pi_{X_+}\widehat D_f.
\]
The first expansion bounds its output perturbation by $C_m\|e\|_2$;
composing the two uniform expansions is therefore legitimate.
Since $\widehat D_fY_f=X_+$, the middle projection is annihilated by
the final projection. Finally $r'p'/(pr)=PQ=d$, proving
\eqref{eq:root-full-derivative} and the uniform expansion.
\end{proof}

\begin{proposition}[Hilbert-space gluing and the single-profile spectrum]\label{prop:root-gluing}
At $f=F$, let $\widehat A_*$ be the first kernel in
\eqref{eq:root-kernel-legs} and put
$B_0=\widehat A_*^*\widehat A_*$.
This operator is compact, self-adjoint, and nonnegative, with
\[
 B_0X=X,\qquad B_0X'=d^{-1}X'.
\]
The scale eigenvalue $d^{-1}$ is simple and all remaining eigenvalues,
apart from the constant-mass eigenvalue $1$, are strictly less than
$d^{-1}$. The principal part of $\mathcal N$ at $X$ is
\begin{equation}\label{eq:root-single-derivative}
 \mathcal L_0=d(B_0-X\otimes X).
\end{equation}
In particular, $\mathcal L_0$ is compact, self-adjoint, and nonnegative;
its eigenvalue $1$ is simple with eigenvector $X'$, and
$\|\mathcal L_0|_{\{X'\}^{\perp}}\|<1$.

For $m\ge2$, write $X_i=S_{t_i}X$, $B_i=S_{t_i}B_0S_{-t_i}$, and
\[
 M_{\mathbf t}
 =d\sum_{i=1}^mB_i
   -\frac d m\Big(\sum_iX_i\Big)\otimes\Big(\sum_iX_i\Big).
\]
Then
\begin{equation}\label{eq:root-operator-gluing}
 \|\mathcal L_{\mathbf t}-M_{\mathbf t}\|_{2\to2}\longrightarrow0
 \qquad\text{as }D(\mathbf t)\longrightarrow\infty.
\end{equation}
\end{proposition}

\begin{proof}
At a single extremizer the two conditional operators are adjoints with
respect to their probability measures, so
$\widehat D_*=\widehat A_*^*$.
Proposition~\ref{prop:direct-tail-hs} gives the Hilbert--Schmidt property
of these kernels. Conjugation by $h\mapsto Xh$ identifies $B_0$ with
the conditional operator in Theorem~\ref{thm:radialspectrum}.
That theorem gives the stated eigenvalues; its scale eigenfunction
$F'/F$ becomes a nonzero multiple of $X'$.
Formula~\eqref{eq:root-full-derivative} then gives
\eqref{eq:root-single-derivative}.

For separated centers, compact-support approximation in the profile
Lebesgue spaces, followed by continuity of the power and normalization
maps, gives
\begin{equation}\label{eq:root-reference-decoupling}
 g-m^{-1/r}\sum_iS_{t_i}H\to0\text{ in }L^r,\qquad
 f_+-m^{-1/p}\sum_iS_{t_i}F\to0\text{ in }L^p,\qquad
 c_f,b_f\to L_m.
\end{equation}
Convolution also gives uniform convergence on every fixed translated
core. There the reference profiles are bounded below, so the kernels
in \eqref{eq:root-kernel-legs} converge uniformly to their single-profile
counterparts. Their amplitude factor
$m^{-1/p+1/2}m^{-1/r+1/2}=m^{-\alpha/n}$ cancels the factor in $L_m$.

To control the complements of the cores, put
$\mathsf U=\sum_iS_{t_i}\mathsf u$ and
$\mathsf V=\sum_iS_{t_i}\mathsf v$.
Finite-sum comparison and $f_+\asymp_m f$ bound the squares of both
kernels by
\[
 C_m\mathcal K(t-z)^2
       \mathsf U(t)^{P-1}\mathsf V(z)^{Q-1}.
\]
If $Q\le1$, this is at most
\[
 C_m\sum_i\mathcal K(t-z)^2
       \mathsf u_i(t)^{P-1}\mathsf v_i(z)^{Q-1}.
\]
Each summand is a common translation of the integrable function in
Proposition~\ref{prop:direct-tail-hs}. Its integral off the paired core
$|t-t_i|,|z-t_i|\le R$ tends to zero uniformly as $R\to\infty$.
If $Q>1$, use instead the sum over all pairs $i,j$.
The diagonal terms have the same bound, while each off-diagonal integral
is a translate of
\[
 \iint\mathcal K(t-z+t_i-t_j)^2
       \mathsf u(t)^{P-1}\mathsf v(z)^{Q-1}\dd t\dd z\longrightarrow0.
\]
Here the kernel is bounded and tends to zero at both ends, and both
profile powers are integrable. Dominated convergence applies.
Consequently,
\[
 \left\|\widehat A_f-\sum_iS_{t_i}\widehat A_*S_{-t_i}\right\|_{\rm HS}
 +\left\|\widehat D_f-\sum_iS_{t_i}\widehat D_*S_{-t_i}\right\|_{\rm HS}
 \longrightarrow0.
\]
Finite-rank approximation and weak vanishing of translations imply
$\|\widehat D_*S_h\widehat A_*\|_{2\to2}\to0$ as $|h|\to\infty$.
Thus
$\|\widehat D_f\widehat A_f-\sum_iB_i\|_{2\to2}\to0$.
Moreover, \eqref{eq:root-reference-decoupling} gives
$X_+-m^{-1/2}\sum_iX_i\to0$ in $L^2$.
Since $B_iX_i=X_i$ and $\|B_iX_j\|_2\to0$ for $i\ne j$,
substitution into \eqref{eq:root-full-derivative} proves
\eqref{eq:root-operator-gluing}.
\end{proof}

\begin{lemma}[Phase choice and control by the update residual]\label{lem:root-residual}
Fix $m\ge2$. For all sufficiently separated centers and sufficiently
close nonnegative unit inputs $x$, the centers can be changed by at most
$C_m\|x-U_{\mathbf t}\|_2$ so that, with the corrected centers,
\[
 e=x-U_{\mathbf t},\qquad
 \langle e,S_{t_i}X'\rangle=0\quad(1\le i\le m).
\]
At these centers one has
\begin{equation}\label{eq:root-residual-bound}
 \|x-U_{\mathbf t}\|_2
 \le C_m\bigl(\|\mathcal N(x)-x\|_2+\eta_{\mathbf t}\bigr),\qquad
 \eta_{\mathbf t}=\|\mathcal N(U_{\mathbf t})-U_{\mathbf t}\|_2.
\end{equation}
All constants are uniform. In particular, an $o(1)$ approximation by
references can be modulated by $o(1)$ corrections; vanishing consecutive
center displacements are preserved.
\end{lemma}

\begin{proof}
First, $X\in H^2(\R)$ and $|X'|+|X''|\le CX$.
Indeed, the kernel logarithmic derivative bound and the stationary
equations give $|F'|\le CF$ and $|H'|\le CH$.
Writing $V=KH$, weak differentiation yields $V''=\mathcal K'*H'$,
so $|V''|\le CV$. Since $F$ is a positive constant times $V^P$,
$|F''|\le CF$, which proves the assertion for $X$.
Also $X'\ne0$, since the nonzero $L^2$ function $X$ cannot be constant.

Writing $V_{\mathbf t}=\sum_iS_{t_i}F$ and
$A_{\mathbf t}=\|V_{\mathbf t}\|_p^{-1}$ gives
$U_{\mathbf t}=(A_{\mathbf t}V_{\mathbf t})^{p/2}$.
The bounds $|F'|+|F''|\le CF$ give uniformly bounded first and second
center derivatives in $L^2$.
Decoupling on fixed cores and the same bounds on their complements give
\[
 \partial_{t_i}U_{\mathbf t}+m^{-1/2}S_{t_i}X'
 \longrightarrow0\quad\text{in }L^2.
\]
Here $A_{\mathbf t}\to m^{-1/p}$ and
$\partial_{t_i}A_{\mathbf t}\to0$, the latter following from
$\int F^{p-1}F'=0$ and separated-profile norm splitting.
For
$\mathscr P_i(\mathbf t,x)=\langle x-U_{\mathbf t},S_{t_i}X'\rangle$,
the center Jacobian at $x=U_{\mathbf t}$ therefore tends to
$m^{-1/2}\|X'\|_2^2 I$.
The derivative has a uniform modulus of continuity, by the regularity
just established. The implicit function theorem, or its contraction
proof with this fixed limiting Jacobian, gives the stated phase choice
and the bound on the correction. Only reference and test functions
are differentiated.

We next prove that $I-M_{\mathbf t}$ is uniformly bounded below on
$Z_{\mathbf t}=\{e:\langle e,S_{t_i}X'\rangle=0\text{ for all }i\}$.
Otherwise there are separating configurations $\mathbf t_j$ and
$e_j\in Z_{\mathbf t_j}$ with $\|e_j\|_2=1$ and
$(I-M_{\mathbf t_j})e_j\to0$.
After a subsequence, let
$S_{-t_{i,j}}e_j\rightharpoonup e_i$ and
$\gamma_i=\langle X,e_i\rangle$.
Compactness of $B_0$ gives
\[
 e_j=d\sum_iS_{t_{i,j}}B_0e_i
   -\frac d m\sum_iS_{t_{i,j}}X\sum_l\gamma_l+o_{L^2}(1).
\]
Translate around center $i$ and pass weakly to the limit:
\[
 e_i=dB_0e_i-\frac d m X\sum_l\gamma_l,\qquad
 \langle e_i,X'\rangle=0.
\]
Taking inner products with $X$ and summing gives
$\sum_l\gamma_l=0$; then $(1-d)\gamma_i=0$, so every $\gamma_i=0$.
The spectral statement in Proposition~\ref{prop:root-gluing} now gives
$e_i=0$: the only eigenvectors of $dB_0$ with eigenvalue $1$ are multiples
of $X'$, which the phase conditions exclude.
The preceding strong representation would therefore imply
$\|e_j\|_2\to0$, a contradiction.
By \eqref{eq:root-operator-gluing}, the same bound holds for
$I-\mathcal L_{\mathbf t}$ on $Z_{\mathbf t}$ once the separation is large.

Finally, for the modulated error, \eqref{eq:root-full-expansion} gives
\[
 \mathcal N(x)-x
 =\mathcal N(U_{\mathbf t})-U_{\mathbf t}
       +(\mathcal L_{\mathbf t}-I)e+o_{L^2}(\|e\|_2).
\]
Apply the lower bound and absorb the uniform remainder to obtain
\eqref{eq:root-residual-bound}.
The last assertion follows because adding $o(1)$ to each center changes
its consecutive displacement by $o(1)$.
\end{proof}

\subsection{Energy gain and exclusion of multiple profiles}

Retain the preceding notation and put $s=2/p\in(1,2)$.
On $L^2\setminus\{0\}$ define
$\mathcal H(x)=E(|x|^s)$.

\begin{lemma}[Energy gain in root variables]\label{lem:root-energy}
Let $x\ge0$, $\|x\|_2=1$, $f=x^s$, and $f_+=\mathcal T(f)$.
Set $y=\mathcal N(x)$, let $g$ be the normalized best response to $f$, and put
$a=E(f)$ and $b=\|Kg\|_{p'}$. Then
\begin{equation}\label{eq:root-young-gap}
 \frac{p-1}{p}\|y-x\|_2^2
 \le 1-\int_\R x^s y^{2-s}\dd t
 =\frac{b-a}{b}
 \le \frac1p\|y-x\|_2^2.
\end{equation}
Moreover,
\begin{equation}\label{eq:root-gradient}
 \nabla\mathcal H(x)
 =s\bigl[b x^{s-1}y^{2-s}-a x\bigr],\qquad
 \|\nabla\mathcal H(x)\|_2\le C\|\mathcal N(x)-x\|_2.
\end{equation}
The constant $C$ depends only on the fixed parameters. On every region where
$E(f)\ge a_*>0$, one also has
\begin{equation}\label{eq:root-energy-gain}
 E(f_+)-E(f)\ge c\|\mathcal N(x)-x\|_2^2,
 \qquad c=a_*\frac{p-1}{p}.
\end{equation}
\end{lemma}

\begin{proof}
The best-response identities give $Kg=b f_+^{p-1}$ and
$a/b=\int x^s y^{2-s}$. For $\theta=1/p\in(1/2,1)$ and $u,v\ge0$,
\[
 (1-\theta)(u-v)^2
 \le \theta u^2+(1-\theta)v^2-u^{2\theta}v^{2-2\theta}
 \le \theta(u-v)^2.
\]
The lower bound follows from weighted arithmetic--geometric mean with
weights $2\theta-1$ and $2-2\theta$, and the upper bound from
$z^{2\theta}\ge 2\theta z-(2\theta-1)$. Integration proves
\eqref{eq:root-young-gap}. Since $g$ is admissible for the maximization at
$f_+$, $E(f_+)\ge b$, while $b\ge a\ge a_*$; this proves
\eqref{eq:root-energy-gain}.

The map $x\mapsto |x|^s$ is $C^1$ from $L^2$ to $L^p$, so
$\mathcal H$ is $C^1$ away from zero. Differentiation gives the formula in
\eqref{eq:root-gradient}. Its inner product with $x$ is zero, so it is
also the gradient along the unit sphere. Finally,
\[
 |u^{s-1}v^{2-s}-u|\le |v-u|\qquad(u,v\ge0).
\]
Using this inequality, \eqref{eq:root-young-gap}, $b\le C_\sharp$, and
$\|y-x\|_2\le2$ proves the gradient bound.
\end{proof}

\begin{lemma}[Quadratic energy comparison]\label{lem:root-energy-comparison}
For sufficiently separated reference parameters and sufficiently close
nonnegative unit inputs $x$, modulated as in Lemma~\ref{lem:root-residual},
\begin{equation}\label{eq:root-energy-comparison}
 |\mathcal H(x)-\mathcal H(U_{\mathbf t})|
 \le C_m\bigl(\|\mathcal N(x)-x\|_2^2+\eta_{\mathbf t}^2\bigr).
\end{equation}
\end{lemma}

\begin{proof}
Write $U=U_{\mathbf t}$, $d=\|x-U\|_2$, and join $U$ to $x$ by the
normalized positive chord
\[
 z(v)=\frac{(1-v)U+vx}{\|(1-v)U+vx\|_2},\qquad 0\le v\le1.
\]
For small $d$, one has $\|z(v)-U\|_2\le Cvd$ and
$\|z'(v)\|_2\le Cd$. The uniform expansion in
Lemma~\ref{lem:root-expansion}, applied at $U$, gives
\[
 \|\mathcal N(z(v))-z(v)\|_2\le\eta_{\mathbf t}+C_mvd.
\]
Integrating \eqref{eq:root-gradient} along this chord yields
\[
 |\mathcal H(x)-\mathcal H(U)|
 \le C_m(\eta_{\mathbf t}d+d^2).
\]
Lemma~\ref{lem:root-residual} bounds
$d\le C_m(\|\mathcal N(x)-x\|_2+\eta_{\mathbf t})$ and therefore gives
\eqref{eq:root-energy-comparison}.
\end{proof}

\begin{proposition}[Exclusion by energy gain]\label{prop:root-exclusion}
Suppose $x_{k+1}=\mathcal N(x_k)$ and
$\mathcal H(x_k)\uparrow L>0$. There cannot be modulated references
$U_k=U_{\mathbf t_k}$, with $D(\mathbf t_k)\to\infty$ and
$\|x_k-U_k\|_2\to0$, and numbers $\varepsilon_k>0$ such that
\begin{equation}\label{eq:root-exclusion-hypotheses}
 \mathcal H(U_k)\ge L+c_0\varepsilon_k,\qquad
 \eta_{\mathbf t_k}^2=o(\varepsilon_k),\qquad
 \frac{\varepsilon_{k+1}}{\varepsilon_k}\longrightarrow1,
\end{equation}
where $c_0>0$ is fixed.
\end{proposition}

\begin{proof}
Put $d_k=L-\mathcal H(x_k)\ge0$ and
$R_k=\|\mathcal N(x_k)-x_k\|_2$. By
Lemma~\ref{lem:root-energy-comparison} and
\eqref{eq:root-exclusion-hypotheses}, for large $k$,
\[
 d_k+\frac{c_0}{2}\varepsilon_k\le C_mR_k^2.
\]
Since $\mathcal H(x_k)\to L>0$, the energy gain
\eqref{eq:root-energy-gain} gives a fixed $\gamma\in(0,1)$ such that
\begin{equation}\label{eq:root-deficit-recurrence}
 d_{k+1}\le(1-\gamma)d_k-\gamma\varepsilon_k.
\end{equation}
Set $w_k=d_k/\varepsilon_k$. Slow variation of $\varepsilon_k$ implies
that, eventually, $w_{k+1}\le\rho w_k-\beta$ for some fixed
$0<\rho<1$ and $\beta>0$. Iteration gives
\[
 0\le w_{k+j}
 \le\rho^j w_k-\frac{\beta(1-\rho^j)}{1-\rho},
\]
which is impossible for large $j$.
\end{proof}

\begin{theorem}[A single concentration profile]\label{thm:it-nomultibubble}
For every orbit satisfying Theorem~\ref{thm:mainiteration}, the number of profiles in
Theorem~\ref{thm:it-profiles} is $m=1$.
\end{theorem}

\begin{proof}
Suppose $m\ge2$. Corollary~\ref{cor:it-centers} gives separated centers
approximating the entire orbit, with consecutive displacements tending
to zero. Apply the modulation of Lemma~\ref{lem:root-residual}; its
$o(1)$ adjustment preserves these properties. Set
\[
 x_k=F_k^{p/2},\qquad U_k=F_{\mathbf t_k}^{p/2},\qquad
 \varepsilon_k=\varepsilon(\mathbf t_k).
\]
Then $\|x_k-U_k\|_2\to0$ and
$\mathcal H(x_k)=E(F_k)\uparrow L_m=C_\sharp m^{-\alpha/n}$.
Proposition~\ref{prop:it-interaction} gives
$\mathcal H(U_k)\ge L_m+c_m\varepsilon_k$ for large $k$, with $c_m>0$.
Moreover, the elementary power inequality and
Proposition~\ref{prop:it-forcing} give
\[
 \eta_{\mathbf t_k}^2
 \le\|\mathcal T(F_{\mathbf t_k})-F_{\mathbf t_k}\|_p^p
 =\zeta(\mathbf t_k)^p=o(\varepsilon_k).
\]
The vanishing center displacements and
\eqref{eq:it-loglipschitzinteraction} imply
$\varepsilon_{k+1}/\varepsilon_k\to1$.
Proposition~\ref{prop:root-exclusion} now gives a contradiction.
If the chosen component is $H_k$, the same argument applies to its full
two-step orbit.
\end{proof}

\begin{proof}[Proof of Theorem~\ref{thm:it-orbitcompact}]
Take any $k_j\to\infty$.
Theorem~\ref{thm:it-profiles} gives a finite-profile representation along this time subsequence; its number $m$ is determined by the common energy limit $L$ of the entire orbit.
Theorem~\ref{thm:it-nomultibubble} gives $m=1$. Hence, after passing to a subsequence, there are $c_j\in\R$ such that
\[
 \|F_{k_j}-S_{c_j}F\|_p+\|H_{k_j}-S_{c_j}H\|_r\longrightarrow0.
\]
Applying $S_{-c_j}$ to both components and using its isometry gives \eqref{eq:it-orbitcompact}, with limit $(F,H)$ and both norms equal to one.
The finite representation, error estimates, and profile count used here have all been derived from the actual recurrence, without assuming convergence of the function sequence beforehand.
\end{proof}

\section{Convergence at a fixed scale}\label{sec:iteration-globalproof}

\begin{proposition}[Sharp energy limit and convergence of the entire shape]\label{prop:it-limitoptimal}\label{prop:it-wholeshape}
Under the assumptions of Theorem~\ref{thm:mainiteration}, one has $L=C_\sharp$ and
\begin{equation}\label{eq:it-wholeshape}
 \inf_{c\in\R}
 \bigl(\|F_k-S_cF\|_p+\|H_k-S_cH\|_r\bigr)\longrightarrow0.
\end{equation}
\end{proposition}
\begin{proof}
Theorem~\ref{thm:it-nomultibubble} gives $m=1$. The quantization formula
\eqref{eq:it-quantization} then gives $L=C_\sharp$, while
Corollary~\ref{cor:it-centers} gives \eqref{eq:it-wholeshape}.
\end{proof}

Shape convergence allows a different translation at each step.
The next argument uses the spectral gap transverse to $F'$ to prove
geometric decay of the error and summability of the translation corrections.
As before, select the full two-step component with $P>1$.

\begin{theorem}[Local attraction to an extremizer at a fixed scale]\label{thm:it-local}
There exists $\delta_*>0$ such that, if a nonnegative unit-norm input
$h\in L^p$ satisfies $\inf_c\|h-S_cF\|_p<\delta_*$, then
$\mathcal T^kh$ converges to a fixed translate $S_{c_\infty}F$,
with $c_\infty\in\R$, and the corresponding half-steps converge to
$S_{c_\infty}H$. Both distances decay geometrically.
\end{theorem}

\begin{proof}
As above, take the full-step component with $1<p<2$.
Put $X=F^{p/2}$ and $x=h^{p/2}$. The inequality
$\|x-S_cX\|_2^2\le\|h-S_cF\|_p^p$ reduces the problem to the
root iteration $\mathcal N$ near $X$; translation equivariance lets us
take the initial translate to be zero.
The single-profile expansion in Lemma~\ref{lem:root-expansion} gives
\[
 \mathcal N(X+e)=X+\mathcal L_0e+o_{L^2}(\|e\|_2).
\]
By Proposition~\ref{prop:root-gluing} and
Theorem~\ref{thm:radialspectrum}, $\mathcal L_0$ is compact,
self-adjoint, and nonnegative. Its eigenvalue $1$ is simple with
eigenvector $X'$, and its norm on $\{X'\}^{\perp}$ is strictly less
than one. Normalization removes the constant-mass direction.

The regularity calculation in Lemma~\ref{lem:root-residual} gives
$X\in H^2(\R)$ and $X'\ne0$.

Define
\[
 \ell(e)=\frac{\langle X',e\rangle}{\|X'\|_2^2},\qquad
 \mathscr P(c,z)=\ell(S_{-c}z-X).
\]
Transferring translation to the test function $X'$ shows that
$\mathscr P$ is $C^1$ on $\R\times L^2$, and
$\partial_c\mathscr P(0,X)=1$.
The implicit function theorem gives a phase $c(z)$ near $X$ satisfying
\[
 \mathscr P(c(z),z)=0,\qquad
 c(z)=-\ell(z-X)+o(\|z-X\|_2),\qquad
 |c(z)|\le C\|z-X\|_2.
\]
Thus $e=S_{-c(z)}z-X$ is orthogonal to $X'$ and
$\|e\|_2\le C\|z-X\|_2$.

For such an error, set $y=\mathcal N(X+e)$ and $\delta=c(y)$.
Since $\ell\mathcal L_0e=0$, we have $\delta=o(\|e\|_2)$ and
\[
 S_{-\delta}y-X=\mathcal L_0e+o_{L^2}(\|e\|_2).
\]
Here $\|S_{-\delta}X-X\|_2=O(|\delta|)$, and
$\|(S_{-\delta}-I)\mathcal L_0\|_{2\to2}\to0$ by compactness and
strong continuity of translations. Choose
$\|\mathcal L_0|_{\{X'\}^{\perp}}\|<\theta<1$.
For every sufficiently small admissible error,
\begin{equation}\label{eq:it-localphase-onestep}
 \|S_{-c(y)}y-X\|_2\le\theta\|e\|_2,\qquad
 |c(y)|\le C\|e\|_2.
\end{equation}

Starting with $c_0=c(x)$ and $e_0=S_{-c_0}x-X$, define
\begin{equation}\label{eq:it-localphase-recursion}
 \begin{aligned}
 y_k&=\mathcal N(X+e_k),&\qquad \delta_k&=c(y_k),\\
 c_{k+1}&=c_k+\delta_k,& e_{k+1}&=S_{-\delta_k}y_k-X.
 \end{aligned}
\end{equation}
Nonnegativity and unit norm are preserved, and the contraction keeps
each error in the same phase neighborhood. Translation equivariance
gives $\mathcal N^kx=S_{c_k}(X+e_k)$, while
\begin{equation}\label{eq:it-localphase-bounds}
 \|e_k\|_2\le\theta^k\|e_0\|_2,\qquad
 |\delta_k|\le C\theta^k\|e_0\|_2.
\end{equation}
Consequently $c_k\to c_\infty\in\R$ and
$|c_k-c_\infty|\le C\theta^k\|e_0\|_2$.
Since $X'\in L^2$,
\begin{equation}\label{eq:it-fixedphase-rate}
 \|\mathcal N^kx-S_{c_\infty}X\|_2
 \le C\theta^k\|e_0\|_2.
\end{equation}
The power map $z\mapsto z^{2/p}$ is locally Lipschitz from
$L^2_+$ to $L^p_+$ because $p<2$, giving the stated geometric
convergence of $\mathcal T^kh$.
Finally the normalized half-step response is locally H\"older from
$L^p_+$ to $L^r_+$ with exponent $\min(1,Q)$, by the HLS inequality
and the scalar estimates for positive powers. It commutes with
translation and sends $F$ to $H$, so its iterates converge geometrically
to the same fixed translate of $H$.
\end{proof}

\begin{proof}[Proof of Theorem~\ref{thm:mainiteration}]
Lemma~\ref{lem:iteration-monotone} proves that the iteration is well defined
and the objective is nondecreasing. Proposition~\ref{prop:it-wholeshape}
identifies its limit as $C_\sharp$ and gives shape convergence.
Thus, at some finite index, the full two-step component selected with
outer power greater than one enters the neighborhood of
Theorem~\ref{thm:it-local}.
The theorem yields a finite limiting translation $c_\infty$ for this
component and its adjacent half-step. If the selected component is $H_k$,
these two limits concern $H_k$ and $F_{k+1}$; shifting the index gives the
same limit for $F_k$.
In either case,
$F_k\to S_{c_\infty}F$ in $L^p$ and
$H_k\to S_{c_\infty}H$ in $L^r$.
The cylindrical isometry, with $\rho_\infty=e^{-c_\infty}$, gives
\eqref{eq:mainiterationlimit} for center zero.
Translation equivariance in Euclidean space restores the initial center $x_0$.
\end{proof}

\end{document}